\documentclass[11pt]{article}
\usepackage{amsmath,amsfonts,amsthm,amscd}

\theoremstyle{plain}
\newtheorem{thm}{Theorem}[section]
\newtheorem{pro}[thm]{Proposition}
\newtheorem{cor}[thm]{Corollary}
\newtheorem{lem}[thm]{Lemma}

\theoremstyle{definition}
\newtheorem{ex}[thm]{Example}
\newtheorem{rem}[thm]{Remark}
\newtheorem{defn}[thm]{Definition}

\newcommand{\Aut}{\operatorname{Aut}}

\newcommand{\End}{\operatorname{End}}
\newcommand{\Hom}{\operatorname{Hom}}

\newcommand{\FF}{{\mathbb{F}}}

\newcommand{\ZZ}{{\mathbb{Z}}}

\newcommand{\Ker}{\mbox{\rm{Ker}}}

\newcommand{\Res}{\mbox{\rm{Res}}}

\newcommand{\gl}{{\mathfrak{gl}}}

\newcommand{\MA}{{\mathbb{A}}}
\newcommand{\MB}{{\mathbb{B}}}
\newcommand{\MI}{{\mathbb{I}}}
\newcommand{\MX}{{\mathbb{X}}}
\newcommand{\MQ}{{\mathbb{Q}}}
\newcommand{\MU}{{\mathbb{U}}}

\def\div{\raise 1pt \hbox{\big|}}

\newcommand{\la}{\langle}
\newcommand{\siff}{\Leftrightarrow}
\newcommand{\ra}{\rangle}
\newcommand{\iso}{\cong}

\def\maprt#1{\smash{\,\mathop{\longrightarrow}\limits^{#1}\,}}

\begin{document}

\title{Quadratic Maps and Bockstein Closed Group Extensions}
\author{Jonathan Pakianathan and Erg\"un
Yal\c c\i n\thanks{Partially supported by a grant
from Turkish Academy of Sciences (T\" UBA-GEB\. IP/2005-16).} }
\maketitle

\begin{abstract}
Let $E$ be a central extension of the form $0 \to V \to G \to W
\to 0$ where $V$ and $W$ are elementary abelian $2$-groups.
Associated to $E$ there is a quadratic map $Q: W \to V$ given by
the $2$-power map which uniquely determines the extension. This
quadratic map also determines the extension class $q$ of the
extension in $H^2(W,V)$ and an ideal $I(q)$ in $H^2(G, \ZZ /2)$
which is generated by the components of $q$. We say $E$ is
Bockstein closed if $I(q)$ is an ideal closed under the Bockstein
operator.

We find a direct condition on the quadratic map $Q$ that
characterizes when the extension is Bockstein closed. Using this
characterization, we show for example that quadratic maps induced
from the fundamental quadratic map $Q_{\gl _n}: \gl _n (\FF _2)
\to \gl _n (\FF _2)$ given by $Q(\MA)= \MA +\MA ^2$ yield
Bockstein closed extensions.

On the other hand, it is well known that an extension is Bockstein
closed if and only if it lifts to an extension $0 \to M \to
\widetilde{G} \to W \to 0$ for some $\ZZ /4[W]$-lattice $M$. In
this situation, one may write $\beta (q)=Lq$ for a ``binding
matrix'' $L$ with entries in $H^1(W, \ZZ /2)$. We find a direct
way to calculate the module structure of $M$ in terms of $L$.
Using this, we study extensions where the lattice $M$ is
diagonalizable/triangulable and find interesting equivalent
conditions to these properties.


\end{abstract}

\section{Introduction}
\label{sect:Introduction} It is of course well known that in
studying the cohomology of finite groups the cohomology of
$p$-groups plays a fundamental part. Any $p$-group $P$ has a
central series with elementary abelian factors (for example the
central Frattini series) and so in principle, its cohomology can
be studied by iteratively considering central extensions with
elementary abelian kernel.

The simplest nontrivial situation is given by $p$-groups $P$ that
fit in a central short exact sequence
$$ 0 \to V \to P \to W \to 0 $$
where $V$ and $W$ are elementary abelian $p$-groups (finite $\FF
_p$-vector spaces).

For $p$ odd, such extensions are in bijective correspondence with a
choice of $p$-power map $\phi: W \to V$ which is linear and an
alternating bilinear commutator map $$[\cdot, \cdot ]: W \wedge W
\to V.$$ For $p=2$, such extensions are in bijective correspondence
with a quadratic $2$-power map $Q: W \to V$ which has an associated
bilinear map given by the commutator. Since the power map is not
linear but quadratic in this case, this adds considerable difficulty
in this situation.

In this paper we study various algebraic properties of these
extensions in the case $p=2$. We start with a review of the above
mentioned facts and an explicit description of how to recover the
extension class in $H^2(W,V)$ of the extension from the quadratic
map $Q$. The components of the extension class form a quadratic
ideal $I(Q) \subseteq H^*(W,\FF _2)$.

In order to compute the cohomology of the group $G$ given by the
extension, an important condition that is often considered is
whether the ideal $I(Q)$ is closed under the Bockstein operator (and
hence under the Steenrod algebra). Our first theorem  gives a
necessary and sufficient condition for the ideal $I(Q)$ to be
Bockstein closed based on conditions on the associated quadratic
form $Q$ of the extension:

\begin{thm}
\label{thm:Intro B-Ptheorem} Let $Q: W \to V$ be a quadratic map, and let
$B$ be the bilinear map associated to $Q$. Then, $Q$ is Bockstein
closed if and only if there exists a bilinear map $P: V \times W \to
V$ which satisfies the identity
$$P(Q(x),y)=B(x,y) +P(B(x,y),x)$$
for all $x,y \in W$.
\end{thm}

While this condition might seem mysterious at first, the case of
$P=B$ happens when $Q$ is obtained from the 2-power map of a
$2$-restricted Lie algebra (see Section
\ref{sect:StronglyB-closed}). A fundamental example of a quadratic
map is $Q_{\gl _n }: \gl _n (\FF_2) \to \gl _n ( \FF _2)$ given by
$$ Q_{\gl _n}(\MA) = \MA + \MA ^2.$$
Using the criterion above, it follows that any quadratic map
induced from $Q_{\gl _n}$ (by restriction of domain and codomain)
has an associated extension which is Bockstein closed.

Let $G$ be the kernel of the mod $2$ reduction map $Gl_n (\ZZ /8)
\to GL _n (\ZZ /2 )$. It is easy to see that $G$ fits into a central
short exact sequence
$$0 \to \gl _n (\FF _2) \to G \to \gl _n (\FF _2) \to 0 $$
with associated quadratic form $Q_{\gl _n}$. It follows that this
extension and its restrictions to suitable subspaces such as
$\mathfrak{sl} _n (\FF _2)$ or $\mathfrak{u} _n (\FF _2)$ are
Bockstein closed. (This is an important ingredient in calculating
the cohomology of these groups as rings and as modules over the
Steenrod algebra. See for example \cite{BrowPak} and
\cite{MinhSym}.)

In the case that $E$ is a Bockstein closed extension, there is a
cohomology class $L$ in $H^1 (W, \End (V))$ such that $\beta
(q)=Lq$ holds. Here the multiplication $Lq$ is given by the
composite
$$ H^1 (W, \End (V))\otimes H^2(W,V) \maprt{\cup} H^3 (W,
\End(V)\otimes V ) \to H^3(W, V)$$ where the second map is induced
by the evaluation map $ev: \End (V)\otimes V \to V$. In this case
we call the element $L$ a {\bf binding operator}. Note that if we
choose a basis for $V$, then we can write $L$ as a matrix with
entries in $H^1 (W, \FF_2)$ and express the equation $\beta
(q)=Lq$ as a matrix equation, where $q$ in the equation denotes
the column matrix whose entries are the components of $q$.

It is well known that the central extension $E: 0 \to V \to G \to
W \to 0$ has a ``uniform'' lift to an extension $\widetilde{E}: 0
\to M \to \widetilde{G} \to W \to 0$, where $M$ is a free $\ZZ
/4$-module with $M/2M = V$ and $\widetilde{G}/2M=G$, if and only
if the extension $E$ is Bockstein closed. (We will reprove this in
this paper when proving more general results.)  The lifting module
$M$ has a $W$-action that is in general nontrivial and this yields
a map $\rho _M : W \to GL(M)$ where $GL(M)$ denotes the
automorphism group of $M$. Note that $\rho _M$ sits in the
subgroup $GL(M;V,V)$ of $GL(M)$ formed by automorphisms of $M$
that induce trivial action on $M/2M$ and $2M$. In the paper we
derive a certain explicit exponent-log correspondence between
$\Hom (W, GL(M; V, V))$ and $H^1(W, \End (V))$, and conclude the
following:

\begin{thm}
\label{thm:Intro mainresult} Let $E: 0 \to V \to G \to W \to 0$ be
a central extension. If $E$ lifts to an extension $\widetilde{E}$
with kernel $M$, then $\log ^*(\rho _M)\in H^1 (W, \End (V))$ is a
binding operator for $E$. Conversely, if  $L$ is a binding
operator for $E$, then $E$ lifts to an extension $\widetilde{E}$
with kernel $M$ where $M$ has a representation $\exp ^* (L)\in
\Hom (W, GL(M; V,V))$.
\end{thm}

This allows one to compute a binding operator $L$ from the module
structure of a lifting module $M$ and vice-versa. Note that in
general the lifting module is not unique, but there is a unique
lifting module for every chosen binding operator. In the paper we
give examples of extensions which have more than one lifting
module and more than one binding operator (see Example
\ref{ex:notunique}).

We then go on in Section~\ref{sect:diagonalizability} to study the
conditions for a Bockstein closed extension to have a diagonalizable
binding operator $L$ (or equivalently a diagonalizable lifting
module $M$). We get the following theorem characterizing the
diagonalizable situation:

\begin{thm}
\label{thm:intro summary} Let $Q : W \to V$ be a Bockstein closed
quadratic map, and let $E: 0 \to V \to G \to W \to 0 $ be the
central extension associated to $Q$. Let $q$ denote the extension
class for $E$. Then, the following are equivalent: \\
(i) $Q$ is diagonalizable, i.e., there is a basis for $V$ such that
the components of $q$ are individually
Bockstein closed. \\
(ii) There is a choice of basis of $V$, such that the components
$q_1, \dots, q_n$ of $q$ all decompose as $q_i=u_iv_i$
where $u_i,v_i$ are linear polynomials. \\
(iii) There exists a diagonalizable $L \in \Hom (W, \End (V))$ such
that $\beta(q)=Lq$.
This is characterized exactly by the equation $\beta(L)+L^2=0$. \\
(iv) $E$ lifts to an extension with kernel $M$ where $M$ is a
direct sum of one dimensional $\ZZ /4 [W]$-lattices. \\
(v) $E$ lifts to an extension with kernel $M$ where
$M$ is a $\ZZ /8 [W]$-lattice. \\
(vi) $E$ lifts to an extension with kernel $M$ where $M$ is a $\ZZ
[W]$-lattice, i.e., $E$ has a uniform integral lifting.
\end{thm}

Some of the implications above are well known to experts on this
subject. The most interesting implication in the above theorem is
that if $E$ has a uniform integral lifting, then the associated
quadratic map is diagonalizable (see
Proposition~\ref{pro:mainthmforlifting}). In other words, the only
extensions we get from integral extensions are the diagonalizable
ones.

In Section~\ref{sect:triangulable}, we consider triangulability of
extensions and find equivalent conditions to triangulability in
terms of the binding operator $L$ and lifting module $M$. We show
that if the extension is triangulable then it has a binding matrix
$L$ such that $\beta (L)+L^2$ is nilpotent. Using this we give
examples of extensions which are not triangulable.

In Section~\ref{sect:Frattini&Effective}, we show that the
components $q_1, \dots, q_n$ of $q$  form a regular sequence as
long as we make the additional assumptions that $\dim(V)=\dim(W)$
and that $Q$ is effective, i.e., that $Q(w)=0$ is zero only when
$w=0$. This turns out to be another important technical tool in
calculating the cohomology of the groups given by these
extensions.

In Section~\ref{sect:2-powerexact}, we consider Bockstein closed
$2$-power exact extensions (see Definition
\ref{defn:2-powerexact}). We show that for a Bockstein closed
$2$-power exact extension the binding operator $L$, and hence
lifting module $M$, is unique.

In Section~\ref{sect:StronglyB-closed}, we show that any
$2$-restricted Lie algebra $\mathfrak{L}$ gives a Bockstein closed
extension of the form $E: 0 \to \mathfrak{L} \to G \to \mathfrak{L}
\to 0$ and further more that there is a bijection between the
$2$-restricted Lie algebras over $\FF _2$ and the Bockstein closed
extensions with $P=B$ in Theorem~\ref{thm:Intro B-Ptheorem}.

In this paper we leave some questions about triangulability of
Bockstein closed $2$-power exact extensions unanswered. These
questions have equivalent versions in terms of $2$-restricted Lie
algebras which may be easier to answer than the original versions.
The statements of these questions can be found in
Sections~\ref{sect:2-powerexact} and~\ref{sect:StronglyB-closed}.

\section{Preliminaries and Definitions}
\label{sect:Preliminaries}

Throughout this section, let $E$ denote a central extension of the
form $$0 \to V \maprt{i} G \maprt{\pi} W\to 0$$ where $V$ and $W$
are elementary abelian $2$-groups. In this section, we will
explain how the standard theory of group extensions and $2$-groups
applies to $E$, then we will make further definitions.

First we will recall how the cohomology class associated to $E$ is
defined. A {\bf transversal} is a function $s: W \to G$ such that
$\pi \circ s=id$ and $s(0)=1$, where $0 \in W$, $1 \in G$ are the
respective identities. Given a transversal $s$, we define $f : W
\times W \to V$ by the formula $$f(x,y)=i^{-1} \bigl
(s(x+y)^{-1}s(x)s(y)\bigr ).$$ We often identify $V$ with its
image under $i$ and do not write $i^{-1}$ in our formulas. It is
well-known that $f$ satisfies the following normalized cocycle
identities:
$$f(x,0)=f(0,x)=0$$
$$ f(x+y,z)+f(x,y)=f(y,z)+f(x,y+z)$$
A function $f:W\times W \to V$ satisfying these identities is
called  a (normalized) {\bf factor set}, and it is known that
normalized factor sets correspond to $2$-cocycles in the
normalized standard cochain complex $C^* (W,V)$. When a different
transversal is chosen, we get a different factor set, say $f'$,
such that $f-f'$ is a coboundary in $C^2(W,V)$. So, for each
extension $E$, there is a unique cohomology class $q \in H^2 (W,
V)$ which is called the {\bf extension class} of $E$. A standard
result in group extension theory  states that up to a suitable
equivalence relation on group extensions the converse also holds
(see for example, Theorem 3.12 in Brown \cite{Brown}):

\begin{pro}
\label{pro:extclassification} There is a one-to-one correspondence
between equivalence classes of central extensions of the form $E : 0
\to V \to G \to W \to 0$ and cohomology classes in $H^2 (W, V)$
where two extensions $E$ and $E'$ are considered equivalent if there
is a commutative diagram of the following form:
$$
\begin{CD}
E: 0 @>>> V @>>> G @>>> W @>>> 0 \\
@. @| @VV{\varphi}V @| @.\\
E': 0 @>>> V @>>> G' @>>> W @>>> 0.\\
\end{CD}
$$
\end{pro}
Now, we consider the group theoretical properties of $G$. Let $[
\, , \, ]$ denote the commutator $[g,h]=g^{-1}h^{-1}gh$, and $g^h
$ denote the conjugation $h^{-1}gh$. Recall the following
identities which are very easy to verify:
\begin{equation}
\label{eqn:group theory}
\begin{split}
[g,hk] &=[g,k][g,h]^k \\
[gh,k] &=[g,k]^h [h,k]\\
(gh)^2 &=g^2h^2 [h,g]^h
\end{split}
\end{equation}
Applying these to the extension group $G$, we see that the
squaring map induces a map $Q: W \to V$ defined as $Q(x)=(\hat
x)^2$, where $\hat x$ denotes an element in $G$ that lifts $x \in
W$. Similarly, the commutator induces a symmetric bilinear map $B:
W \times W \to V$ defined as $B(x,y)=[\hat x, \hat y]$. Note also
that since $V$ is a central elementary abelian 2-group, the third
row of $(1)$ gives $$B(x,y)=Q(x+y)+Q(x)+Q(y)$$ hence $Q$ is a
quadratic map and $B$ is its associated bilinear map.

We now will explain the relation between the quadratic map $Q$ and
the extension class $q$. Choosing a basis $\{ v_1, \dots, v_n \}$
for $V$, we can write $Q=(Q_1, \dots, Q_n)$ and  $B=(B_1, \dots,
B_n)$. Also using the isomorphism $H^2(W, V) \iso H^2 (W, \ZZ
/2)\otimes V$, we can write $q=(q_1, \dots, q_n)$. We have the
following:

\begin{pro}
\label{pro:coefficients} Let $\{ w_1,\dots , w_m \}$ be a basis for
$W$, and let $\{ x_1,\dots , x_m \}$ denote its dual basis. Under
the isomorphism  $H^* (W, \ZZ /2) \iso \ZZ /2 [x_1, \dots, x_m]$,
the cohomology class $q_k$ is equal to the quadratic polynomial
associated to $Q_k$ for each $k$. In other words, as a polynomial in
$x_i$'s, the cohomology class $q_k$ is of the following form:
$$q_k =\sum _i Q_k (w_i) x_i^2 + \sum _{i<j} B_k (w_i, w_j)
x_i x_j.$$
\end{pro}

\begin{proof} Let $E: 0 \to V \to G \to W \to 0$ be an extension
with quadratic map $Q$ and bilinear map $B$. We define $f: W\times
W \to V$ as the bilinear map such that
\[f(w_i, w_j)=
\begin{cases}
B( w_i , w_j )  & \text{if \: \: $i < j$}, \\
Q(w_i)  &  \text{if \: \: $i=j$}, \\
0 \ \  & \text{if \: \: $i>j$}.
\end{cases}\]
It is clear that, being a bilinear map, $f: W \times W \to V$
satisfies the factor set conditions. Also note that for every $k$,
and for every $w, w' \in W$, we have
$$f_k (w,w') = \sum _i  Q_k(w_i) x_i (w) x_i(w') + \sum_{i<j}  B_k(
w_i, w_j)x_i(w) x_j(w'),$$ so the associated cohomology class for
$f_k$ is equal to
$$\sum _i  Q_k(w_i) x_i ^2 + \sum_{i<j}  B_k(
w_i, w_j)x_i x_j$$ under the isomorphism $H^*(W, \ZZ/2)\iso \ZZ /2
[x_1, \dots, x_m]$.

To complete the proof, we just need to show that $f$ is a factor
set for $E$. Note that for every $w = a_1 w_1 + \cdots + a_m w_m$
in $W$, we have
\begin{equation*}
f(w,w)=\sum  _{i,j} a_i a_j f(w_i, w_j)=\sum
_i a_i Q(w_i) + \sum _{i<j} a_ia_j B(w_i, w_j)=Q(w).
\end{equation*}
Hence, the proof follows from the following lemma.
\end{proof}

\begin{lem}
\label{lem:factorset} Let $E: 0 \to V \to G \to W \to 0$ be an
extension with quadratic map $Q$. A factor set $f:W\times W \to V$
is a factor set for $E$ if and only if $f(w,w)=Q(w)$ for all $w
\in W$.
\end{lem}

\begin{proof} If $f$ is a factor set for $E$ determined by the
splitting $s$, then $$f(w,w)=s(w+w)^{-1}s(w)s(w)=(s(w))^2=Q(w)$$
holds for all $w \in W$. For the converse, let $f$ be a factor set
satisfying $f(w,w)=Q(w)$ for all $w \in W$. Suppose $f'$ is a
factor set for $E$. Then, $f''=f+f'$ is a factor set such that
$f''(w,w)=0$ for all $w\in W$. If $0 \to V \to G'' \maprt{\pi} W
\to 0$ is the extension with factor set $f''$ and transversal $t:W
\to G''$, then for every $g \in G''$, we have $g^2=[t(\pi (g))
]^2=f''(\pi(g), \pi(g))=0$. Thus $G''$ is of exponent $2$, and
hence an elementary abelian $2$-group. This implies that the
extension $0 \to V \to G'' \maprt{\pi} W \to 0$ splits, and
therefore $f''$ is cohomologous to zero. This shows that $f'$ is
also a factor set for $E$.
\end{proof}

As a consequence of Propositions~\ref{pro:extclassification}
and~\ref{pro:coefficients}, we obtain the following:

\begin{cor}
\label{cor:conclusion} Given a quadratic map $Q: W \to V$, there
is a unique (up to equivalence) central extension
$$E(Q): 0 \to V \to G(Q) \to W \to 0$$ with $Q(w)=(\hat w )^2$ for
all $w \in W$.
\end{cor}

In fact, there is a natural equivalence between the category of
quadratic maps $Q: W \to V$ and the category of central extensions
of the form $E: 0 \to V \to G\to W \to 0.$ We skip the details of
this categorical equivalence, since we do not need this for this
paper.

We end this section with some examples:

\begin{ex}
\label{ex:gl} Let $\gl _n  (\FF _2 )$ denote the vector space of
$n \times n$ matrices over $\FF _2$. Define $Q: \gl _n (\FF _2)
\to \gl _n (\FF _2) $ by $Q(\MA)=\MA + \MA ^2$. Then one computes
$$Q(\MA +\MB )=Q(\MA )+Q(\MB ) +[\MA ,\MB ]$$ where $[\MA ,\MB
]=\MA \MB + \MB \MA$. Thus $Q$ is a quadratic map with associated
bilinear map equal to the Lie bracket $[\cdot , \cdot ]$. We
denote this quadratic map by $Q_{\gl _n}$.

Let $K_n (\ZZ /2^m)$ be the kernel of the map $GL_n(\ZZ /2^m) \to
GL_n(\ZZ /2)$ defined by the mod $2$ reduction of the entries of
the given matrix. Observe that
$$K_n (\ZZ /8)=\{ \MI+2 \MA \ | \ A \in \gl _n (\ZZ /4) \}$$
is a non-abelian
group of exponent $4$. Consider the mod $4$ reduction map
$$\varphi : K _n (\ZZ / 8 ) \to K _n (\ZZ /4).
$$ Note that both
$ K_n(\ZZ /4)$ and $\Ker ( \varphi )$ are isomorphic to $\gl _n
(\FF _2 )$, and that $\ker \varphi$ is a central subgroup of $K_n
(\ZZ /8 )$. So, we get a central extension of the form
$$ E : 0 \to  \gl _n (\FF _2)  \to K_n (\ZZ /8) \to  \gl _n (\FF _2 )
\to 0.$$ Note that $(\MI +2\MA )^2=\MI +4(\MA +\MA ^2)$, and so the
associated $Q$ for this extension is indeed equal to $Q_{\gl _n}$.
\end{ex}

\begin{ex}
\label{ex:un} If $W$ is a $\FF _2$-subspace of $\gl _n (\FF _2 )$,
such that $\MA \in W $ implies $\MA ^2 \in W$, then $Q_{\gl _n}$
induces a quadratic map $Q: W \to W$. We call such a subspace $W$
square-closed. Since $[\MA ,\MB ]=Q(\MA +\MB ) + Q(\MA) + Q(\MB)$
we have that $W$ is automatically a sub-Lie algebra of $\gl _n (
\FF _2 )$. Some examples of subspaces $W$ with this property are:

(a) $W=\mathfrak{sl} _n (\FF _2 )$, the matrices of trace zero.

(b) $W=\mathfrak{u} _n (\FF _2 )$, the strictly upper triangular
matrices.

\noindent In this case we call the corresponding quadratic map $Q:
W \to W$ a {\bf $\gl$-induced quadratic map}.
\end{ex}

\section{Bockstein Closed Extensions}
\label{sect:B-closed}

We start with some definitions.

\begin{defn}
\label{defn:I(Q)} If $q_1,\dots,q_n$ are components of $q$ with
respect to some basis of $V$ we denote by $I(Q)$ the ideal
$(q_1,\dots,q_n)$ in $H^*(W,\FF _2)$.
\end{defn}

It is easy to see that the ideal $I(Q)$ is indeed independent of the
basis chosen for $V$, and hence is completely determined by $Q$.

\begin{defn}
\label{defn:Bocksteinclosed} We say the quadratic map $Q:W \to V$
is Bockstein closed if $I(Q)$ is invariant under the Bockstein
operator on $H^*(W, \ZZ /2)$. An extension $E(Q) : 0 \to V \to
G(Q) \to W \to 0$ is called Bockstein closed if the associated
quadratic map $Q$ is Bockstein closed.
\end{defn}

Since $I(Q)$ is an ideal generated by homogeneous polynomials of
degree $2$, it is closed under the higher Steenrod operations, so
$I(Q)$ is a Steenrod closed ideal. We will later use this fact in
Section \ref{sect:Frattini&Effective} when we are studying effective
extensions.

The main examples of Bockstein closed extensions are $\gl$-induced
quadratic maps. We will prove the Bockstein closeness of these
extensions at the end of this section. We first start with an easy
observation.

\begin{pro}
\label{pro:L} Let $Q: W \to V$ be a quadratic map, and let $q \in
H^2(W, V)$ be the corresponding extension class. Then, $Q$ is
Bockstein closed if and only if there is a one dimensional class $L
\in H^1 (W, \End (V))$ such that $\beta (q)=Lq$. Here, the
multiplication $Lq$ is given by the composite
$$ H^1 (W, \End (V))\otimes H^2(W,V) \maprt{\cup} H^3 (W,
\End(V)\otimes V ) \to H^3(W, V)$$ where the second map is induced
by the evaluation map $ev: \End (V)\otimes V \to V$, given by
$f \otimes u \to f(u)$.
\end{pro}

\begin{proof} Choose a basis ${\cal B}_W$ for $W$ and ${\cal B}_V$
for $V$, and
let $\{x_1, \dots, x_m \}$ be the dual basis to ${\cal B}_W$. We
can write $q$ as a column vector with entries in $H^2(W, \FF _2)$.
Then, $q$ is Bockstein closed if and only if there is an $n \times
n$ matrix $L$ with entries in linear polynomials in $x_i$'s such
that $\beta (q)=Lq$ where $n=\dim (V)$. Note that $L$ can be
considered as an element in $H^1(W, M_{n\times n} (\FF _2 ))$.
Using the basis ${\cal B}_V$, we can also identify $M_{n\times n}
(\FF _2 )$ with $\End (V)$. So, $L$ can be considered as a class
in $H^1(W, \End(V))$, and it is easy to see that under these
identifications, matrix multiplication $Lq$ corresponds to the map
given above.
\end{proof}

\begin{defn}
\label{defn:bindingoperator} Let $Q: W \to V$ be a quadratic map,
and let $q \in H^2(W, V)$ be the corresponding extension class. If
$L \in H^1 (W, \End (V))$ is a class satisfying $\beta (q)=Lq$,
then we say $L$ is a {\bf binding operator} for $Q$. If a specific
basis for $V$ is chosen then $L$ can represented as a matrix with
linear polynomial entries. In this case we call $L$ a {\bf binding
matrix}.
\end{defn}

It is clear that one can have many different matrices for the same
operator. Any two such matrices will  be conjugate to each other by
a scalar matrix. On the other hand, one can also have two
different binding operators $L_1$ and $L_2$ for a quadratic map
$Q$. The following example shows a case where this happens.

\begin{ex}
\label{ex:notunique} Let $q=(xy, yz)$. Then, $\beta(q)=(xy(x+y),
yz(y+z)).$ We can write
$$\left[\begin{matrix}\beta(q_1)\cr \beta(q_2)\cr\end{matrix}\right]
=\left[\begin{matrix}x+y & 0\cr 0 & y+z \cr\end{matrix}\right]
\left[\begin{matrix}q_1\cr q_2\cr\end{matrix}\right]$$ or
$$\left[\begin{matrix}\beta(q_1)\cr \beta(q_2)\cr\end{matrix}\right]
=\left[\begin{matrix}x+y+z & x\cr z & x+y+z \cr\end{matrix}\right]
\left[\begin{matrix}q_1\cr q_2\cr\end{matrix}\right].$$ Note that
the second binding matrix is not conjugate to a diagonal matrix.
So, there exists more than one binding operator in this case. This
can happen because the ideal $(q_1,q_2)$ is not a free
$k[x,y,z]$-module over the generators $\{q_1,q_2\}$. We will see
later that under stronger conditions there is a unique binding
operator $L$ for $q$.
\end{ex}

The following is a reformulation of Proposition~\ref{pro:L} which
is quite useful in many instances.

\begin{thm}
\label{thm:B-Ptheorem} Let $Q: W \to V$ be a quadratic map, and let
$B$ be the bilinear map associated to $Q$. Then, $Q$ is Bockstein
closed if and only if there exists a bilinear map $P: V \times W \to
V$ which satisfies the identity
$$P(Q(x),y)=B(x,y) +P(B(x,y),x)$$
for all $x,y \in W$.
\end{thm}

Before the proof we first do some calculations. Let $\{ v_1,
\dots, v_n \}$ be a basis for $V$, and let $\{ w_1,\dots , w_m \}$
be a basis for $W$ with dual basis $\{x_1,\dots , x_m \}$. As
usual we write $q=(q_1, \dots, q_n )$ for the extension class $q
\in H^2 (W, V)$ associated to $Q$. Note that by
Proposition~\ref{pro:coefficients}, for each $k=1,\dots, n$, we
have
$$q_k= \sum _i Q_k(i) x_i ^2 + \sum _{i<j} B_k(i, j) x_i
x_j $$ where  $Q_k (i)=Q_k (w_i)$ and $B_k(i,j)=B_k (w_i, w_j)$.
Applying the Bockstein we get
\begin{equation}
\label{eqn:betaqk}
\begin{split}
\beta (q_k) &= \sum _{i<j } B_k(i,j) \bigl ( x_i^2 x_j + x_i x_j
^2 \bigr ) =\sum _{i,j} B_k(i,j) x_i^2 x_j
\end{split}
\end{equation}

Given an element $L \in H^1 (W, \End (V))$, we can consider it as
a homomorphism $L : W \to \End (V) $ via the isomorphism
$$H^1 (W, \End (V)) \iso \Hom (W, \End (V) ).$$ This allows us to
describe a correspondence between classes $L \in H^1 (W, \End
(V))$ and bilinear maps $P : V \times W \to V$. The correspondence
is given by $P(v, w)= L(w) (v)$ for all $v \in V$ and $w \in W$.
Therefore, if $L$ denotes the matrix for an element in $H^1 (W,
\End (V))$ with respect to the above choice of basis for $V$ and
$W$, then we can write
$$L_{ks}=\sum _p P_k (s,p) x_p $$ where $P_k (s,p)$ is short for
$P_k (v_s, w_p)$ and $P_k$ is the $k$-th coordinate of $P$. We
have the following:

\begin{lem}
\label{lem:twoequations} The equation $\beta (q)=Lq$ holds if and
only if the following two equations hold
\begin{equation}
\label{eqn:1}
 P( Q(w_i), w_j) = B (w_i, w_j ) + P ( B(w_i, w_j), w_i)
\end{equation}
\begin{equation}
\label{eqn:2} P(B(w_i,w_j),w_k)+P(B(w_j ,w_k), w_i)+P(B(w_k ,w_i),
w_j)=0
\end{equation}
for all $i,j,k \in \{ 1, \dots, n \}.$
\end{lem}

\begin{proof} First note that by Equation
\ref{eqn:betaqk} above, we have
\begin{equation*}
\beta (q_k)= \sum _{i,j} B_k(i,j) x_i^2 x_j
\end{equation*}
for all $k$. On the other hand, the $k$-th entry of $Lq$ is equal to
$\sum _s L_{ks} q_s $. Writing this sum in detail, we get
\begin{equation*}
\begin{split}
\sum_s L_{ks} q_s &= \sum _s \Bigl [ \Bigl (\sum _p P_k(s,p) x_p
\Bigr ) \Bigl ( \sum_{i} Q_s(i) x_i^2 + \sum_{i<j} B_s(i,j) x_i x_j
\Bigr )
\Bigr ] \\
&= \sum _s  \sum _{p,i} Q_s(i) P_k (s,p) x_i^2 x_p + \sum_s \sum
_{p,\: i<j} B_s (i,j) P_k (s,p) x_p x_i x_j \\
&= \sum _{p,i} \bigl ( \sum_s Q_s(i) P_k (s,p) \bigr ) x_i^2 x_p +
\sum_{p,\: i<j} \bigl ( \sum_s B_s (i,j) P_k (s,p) \bigr )
x_p x_i x_j \\
&= \sum _{i,j} \bigl ( \sum_s Q_s(i) P_k (s,j) \bigr ) x_i^2 x_j +
\sum_{p,\: i<j} \bigl ( \sum_s B_s (i,j) P_k (s,p) \bigr ) x_px_ix_j.\\
\end{split}
\end{equation*}
The second term on the righthand side (SLHS)  can be manipulated
more:
\begin{equation*}
\begin{split}
SLHS &= \sum_{p=i<j} \bigl ( \sum_s B_s (i,j) P_k (s,i) \bigr )
x_i^2 x_j +
\sum_{i<j=p} \bigl ( \sum_s B_s (i,j) P_k (s,j) \bigr ) x_i x_j^2 \\
& + \sum _{\substack{i<j \\ p\neq i,j }}  \bigl ( \sum_s B_s (i,j)
P_k(s,p) \bigr ) x_px_ix_j \\
& = \sum_{i,j} \bigl ( \sum_s B_s (i,j) P_k (s,i) \bigr ) x_i^2
x_j+ \sum _{\substack{i<j \\ p\neq i,j }}  \bigl ( \sum_s B_s
(i,j) P_k(s,p) \bigr ) x_px_ix_j
\end{split}
\end{equation*}
where the last line follows from the identity $B_s (i,j)=B_s
(j,i)$. Putting these equations together, we get
\begin{equation*}
\begin{split}
\beta (q_k)+\sum_s L_{ks} q_s &=  \sum _{i,j} \Bigl ( B_k(i,j)+
\sum_s Q_s(i) P_k (s,j) +
\sum_s B_s (i,j) P_k (s,i)  \Bigr )  x_i^2 x_j  \\
&+ \sum _{\substack{i<j \\ p\neq i,j }}  \bigl ( \sum_s B_s (i,j)
P_k(s,p) \bigr ) x_px_ix_j.\\
\end{split}
\end{equation*}
Fixing an order for $\{p,i,j\}$, the second summand in the above
equation becomes
$$\sum _{p<i<j} \Bigl [ \sum_s \Bigl
( B_s (i,j) P_k(s,p) + B_s (p,j) P_k(s,i)+B_s (p,i) P_k(s,j) \Bigr
) \Bigr ] x_px_ix_j $$ So, we obtain that $\beta (q)+Lq =0$ if and
only if
$$B_k(i,j)+ \sum_s Q_s(i) P_k (s,j) + \sum_s B_s (i,j) P_k
(s,i)=0$$ and
$$\sum _s \Bigl ( B_s (i,j) P_k(s,p) +  B_s (p,j) P_k(s,i)+
B_s (p,i) P_k(s,j) \Bigr ) =0$$ hold for all $i, j, k $. It is
easy to see that these equations are the same as the
Equations~\ref{eqn:1} and~\ref{eqn:2}.
\end{proof}

We will also need the following lemma:

\begin{lem}
\label{lem:lieidentity} If $P: V \times W \to V$ be a bilinear
form satisfying the identity $$P(Q(x),y)=B(x,y) +P(B(x,y),x)$$ for
all $x,y \in W$, then it satisfies
\begin{equation}
\label{eqn:lieidentity}
P(B(x,y),z)+P(B(y,z),x)+P(B(z,x),y)=0
\end{equation}
for all $x,y,z \in W.$
\end{lem}

\begin{proof} We show this by a direct calculation.
\begin{equation*}
\begin{split}
P(B(x,y),z)&= P(Q(x)+Q(y)+Q(x+y),z) \\
& =P(Q(x),z) +P(Q(y),z) + P(Q(x+y),z) \\
&= B(x,z)+P(B(x,z),x)+B(y,z)+P(B(y,z),y)\\
&+ B(x+y, z)+
P(B(x+y,z),x+y) \\
&=P(B(x,z),x) + P(B(y,z), y) + P(B(x+y,z), x+y) \\
&=P(B(x,z),x) + P(B(y,z), y) + P(B(x,z)+B(y,z), x+y)\\
&=P(B(x,z),y) + P(B(y,z),x) \\
&=P(B(y,z),x) + P(B(z,x),y)
\end{split}
\end{equation*}
\end{proof}

Now, we are ready to prove Theorem \ref{thm:B-Ptheorem}.
\begin{proof}[Proof of Theorem~\ref{thm:B-Ptheorem}]
First assume that there is a bilinear map $P: V \times W \to V$
satisfying the identity $P(Q(x),y)=B(x,y) +P(B(x,y),x)$ for all $x,y
\in W$. Let $L$ be the associated cohomology class. By
Lemma~\ref{lem:twoequations}, we just need to show that the bilinear
map $P$ satisfies the Equations~\ref{eqn:1} and~\ref{eqn:2}. For the
first equation take $x=w_i$ and $y=w_j$, and consider the $k$-th
coordinates. To get the second equation, we first use
Lemma~\ref{lem:lieidentity}, and put $x=w_i, y=w_j$ and $z=w_p$ to
the Equation~\ref{eqn:lieidentity}.

For the converse, assume that there is a binding operator $L \in
H^1 (W, \End (V))$ satisfying the equation $\beta (q)=Lq$. Let $P
: V \times W \to V$ be the bilinear map associated to $L$.
Lemma~\ref{lem:twoequations} gives us that $P$ satisfies
Equations~\ref{eqn:1} and~\ref{eqn:2} for all $i,j,k \in \{ 1,
\dots, n \}.$ By linearity this implies that the equation
$$P(Q(x),y)=B(x,y) +P(B(x,y),x)$$ is satisfied for all $x,y \in
W$. Thus the proof is complete.
\end{proof}

As an immediate corollary of Theorem \ref{thm:B-Ptheorem}, we
obtain
\begin{cor}
\label{cor:glinduced} Any $\gl$-induced quadratic map $Q: W \to W$
is Bockstein closed.
\end{cor}
\begin{proof}
For $Q: \gl _n (\FF _2 ) \to \gl _n (\FF _2 )$ we can set $P( \MA
,\MB )=[\MA, \MB]$. Then we compute
\begin{align*}
\begin{split}
P(Q(\MA),\MB)&=[\MA+\MA^2,\MB] \\
&=[\MA,\MB] + [\MA^2,\MB] \\
&=[\MA,\MB] + [[\MA,\MB],\MA]
\end{split}
\end{align*}

Thus equation in Theorem~\ref{thm:B-Ptheorem} holds with $P=B$ and
hence the $Q_{\gl _n}$ is Bockstein closed. It is easy to see $P:
W \times W \to W$ for any square-closed subspace $W$ and so every
$\gl$-induced $Q: W \to W$ is Bockstein closed.
\end{proof}

\section{Uniform Lifting and Binding Operators}
\label{sect:uniformlifting}

Throughout this section $Q:W \to V$ will denote an arbitrary
quadratic form, and $E(Q): 0 \to V \to G(Q) \to W \to 0$ will be
the central extension associated to $Q$. As usual we will denote
the extension class of this extension by $q$.

\begin{defn}
\label{defn:uniformlifting} Let $M$ be a $\ZZ [W]$-module such that
$M/ 2M \iso V$ as $\ZZ /2[W]$-modules. We say $E(Q)$ {\bf lifts to
an extension with kernel $M$} if there is an extension of the form
$$\widetilde E(Q): 0 \to M \to \widetilde G (Q) \to W \to 0$$ such
that the following diagram commutes
$$
\begin{CD}
@. 2M @= 2M @.  \\
@. @VVV @VVV @.  \\
\widetilde E(Q) : \  @. M @>>> \widetilde G (Q)  @>>> W \\
@. @VVV @VVV @| \\
E(Q) : \ @. V @>>> G(Q) @>>> W. \\
\end{CD}
$$
In the case $M$ is a $\ZZ/4$-free $\ZZ /4[W]$-module ($\ZZ
/4[W]$-lattice), we say $E(Q)$ has a {\bf uniform lifting}. If $M$
is a $\ZZ$-free $\ZZ [W]$-module ($\ZZ[W]$-lattice) then we say $E$
has a {\bf uniform integral lifting}.
\end{defn}

It is well known that a quadratic map $Q$ is Bockstein closed if and
only if the associated extension $E(Q)$ has a uniform lifting.
However, the known proofs of this statement do not provide an
explicit connection between the quadratic form and the $\ZZ
/4[W]$-lattice structure of $M$. In this section we introduce the
concept of binding operators and find a direct way to calculate the
module structure of $M$ from the binding operator $L$.

We start with the following observation:

\begin{lem}
\label{lem:longexactsequence}  Suppose that $M$ is a $\ZZ
/4[W]$-lattice such that $2M$ and $M/2M$ are trivial $\ZZ
/2[W]$-modules isomorphic to $V$. Then, $E(Q)$ lifts to an
extension with kernel $M$ if and only if $\delta(q)=0$ where
$\delta : H^2 (W, V) \to H^3 (W, V)$ is the connecting
homomorphism of the long exact sequence for the coefficient
sequence $0 \to V \to M \to V \to 0$ of $\ZZ [W]$-modules.
\end{lem}

\begin{proof} Recall that the long exact sequence for the extension
$0 \to V \to M \to V \to 0 $ is of the form
$$ \cdots \to H ^2 (W, M ) \to H^2 (W, V ) \maprt{\delta} H^3 (W, V
)  \to \cdots $$ So, from this it is clear that $q$ lifts to a
class $\tilde q \in H^2(W, M)$ if and only if $\delta(q)=0$.
Lifting the extension class is equivalent to lifting the
associated extension, so the proof of the lemma is complete.
\end{proof}

So, it remains to understand the boundary homomorphism for the
extension $0 \to V \to M \to V \to 0$. For this we look at the
boundary homomorphism on the chain level. We will consider the
boundary operators associated to sequences of $\ZZ [W]$-modules $0
\to A \to B \to C \to 0$ where $A$ and $C$ have trivial $\ZZ
[W]$-module structure. To do this we will discuss the notions of
binding functions and binding operators.

Let $0 \to A \to B \overset{\pi}{\to} C \to 0$ be a short exact
sequence of abelian groups. Let $GL(B;A,C)$ be the group of
automorphisms of $B$ which induce the identity map on $A$ and $C$.

\begin{lem}
\label{lem:isomorphism} There are isomorphisms $\log: GL(B;A,C) \to
\Hom(C,A)$ and $\exp: \Hom(C,A) \to GL(B;A,C)$ (inverse to each
other) given by $\log(f)(c)=f(\hat{c})-\hat{c}$ for all $c \in C$
where $\hat{c} \in B$ has $\pi(\hat{c})=c$ and
$\exp(\mu)(b)=b+\mu(\pi(b))$ for all $b \in B$.
\end{lem}

\begin{proof} We first show that $\log$ and $\exp$ are well-defined
functions. Note that if $f \in GL(B;A,C)$ then $f(\hat{c})-\hat{c}$
projects trivially under $\pi$ since $f$ induces the identity map on
$C$ and hence $f(\hat{c})-\hat{c} \in A$. Also since $f$ induces the
identity map on $A$, $f(\hat{c})-\hat{c}$ is indeed independent of
the lift $\hat{c}$ of $c$. Finally
$\log(f)(c+d)=f(\hat{c}+\hat{d})-(\hat{c}+\hat{d})=\log(f)(c)+\log(f)(d)$
since $f$ is a homomorphism. Thus $\log(f) \in \Hom(C,A)$.

It is clear that $\exp(\mu)$ defines an endomorphism of $G$ which
induces the identity map on $A$ and $C$. Since $\exp(-\mu)$ is
easily seen to be its inverse, $\exp(\mu)$ is an element of
$GL(B;A,C)$.

Finally since it is easy to see that $\exp$ and $\log$ are inverse
functions, to show they are isomorphisms, we need only check that
$\log$ is a homomorphism. If $f,g \in GL(B;A,C)$ then
\begin{align*}
\begin{split}
\log(f \circ g)(c) &=f(g(\hat{c}))-\hat{c} \\
&=f(g(\hat{c}))-g(\hat{c})+g(\hat{c})-\hat{c} \\
&=\log(f)(c) + \log(g)(c)
\end{split}
\end{align*}
where the final step follows since $g(\hat{c})$ is also a lift of
$c$ since $g$ induces the identity map on $C$. Thus $\log(f \circ
g) = \log(f) + \log(g)$ and the proof is complete.
\end{proof}

Now given a sequence of $\ZZ [W]$-modules $0 \to A \to B \to C \to
0$ such that $A$ and $C$ are trivial $\ZZ [W]$-modules, it is
clear we get a representation $\rho_B: W \to GL(B;A,C)$. Thus we
obtain a homomorphism $\log(\rho_B): W \to \Hom(C,A)$ where
$\log(\rho_B)=\log \circ \rho_B$.

The extension above yields a long exact sequence whose boundary
operator is $\delta: H^*(W, C) \to H^{*+1}(W , A)$. On the other
hand we may also consider the sequence $0 \to A \to B \to C \to 0$
as a sequence of trivial $\ZZ [W]$-modules and this will yield
another long exact sequence with boundary operator $\delta_{triv}:
H^*(W, C) \to H^{*+1}(W, A)$.

In general $\delta$ will be different from $\delta_{triv}$ due to the
``twisting'', i.e., nontrivial action of $W$ on $B$. The next
proposition makes the connection between these two boundary operators
explicit.

\begin{pro}
\label{pro:connectinghomomorphism} Let $0 \to A \to B \to C \to 0$
be a short exact sequence of $W$-modules such that $A$ and $C$ are
trivial $W$-modules. Let $\delta: H^*(W, C) \to H^{*+1}(W, A)$ be
the boundary operator associated to this sequence and
$\delta_{triv}: H^*(W, C) \to H^{*+1}(W, A)$ be the boundary
operator associated with the same sequence and trivial $W$-action.
Then,  $$\delta=\delta_{triv} + \mathbb{M}$$ where $\mathbb{M}$ can
be defined on cochains $\mathbb{M}: C^n(W, C) \to C^{n+1}(W, A)$ as
follows:
$$\mathbb{M}(f)(z_0,\dots,z_n)=\log(\rho_B(z_0))(f(z_1,\dots,z_n))$$
for all $f \in C^n(W, C)$. We will call $\mathbb{M}: H^n(W, C) \to
H^{n+1}(W, A)$ the binding operator of the extension.
\end{pro}

\begin{proof}
Consider normalized standard bar resolutions. The exact sequence
$0 \to A \to B \to C \to 0$ induces a commuting diagram
$$
\begin{CD}
0 @>>> C^n (W,A) @>j>> C^n(W, B ) @>{\pi}>> C^n (W, C) @>>> 0 \\
@. @VV{d }V  @VV{d }V  @VV{d }V @. \\
0 @>>> C^{n+1} (W,A) @>j>> C^{n+1} (W, B ) @>\pi>> C^{n+1} (W, C) @>>>
0 \\
\end{CD}
$$
and connection homomorphism $\delta: C^n(W, C) \to C^{n+1}(W, A)$
defined as follows: For $f \in C^n(W, C)$, we lift (valuewise) to
get $\hat{f} \in C^n(W, B)$. Then $\delta(f)=d\hat{f}$. Using
$\rho_B: W \to GL(B;A,C)$ to explicitly write the $W$ action on
$B$ and computing we get:
\begin{equation*}
\begin{split}
(\delta f)(z_0, \dots , z_n) &= \rho_B(z_0)(\hat{f}(z_1, \dots , z_n))\\
& -\hat{f}(z_0z_1, z_2, \dots , z_n ) \\
&\cdots \\
&\pm \hat{f}(z_0, z_1, \dots, z_{n-1}z_{n} ) \\
&\mp \hat{f}(z_0, z_1, \dots , z_{n-1} )
\end{split}
\end{equation*}
which we write as
\begin{align*}
\begin{split}
(\delta f)(z_0, \dots , z_n)
&= \rho_B(z_0)(\hat{f}(z_1, \dots , z_n)) - \hat{f}(z_1, \dots , z_n)
+ (\delta_{triv} f)(z_0, \dots , z_n) \\
&= \log(\rho_B(z_0))(f(z_1,\dots,z_n)) + (\delta_{triv}
f)(z_0,\dots,z_n) \\
&= (\delta_{triv} f + \mathbb{M}f)(z_0,\dots,z_n)
\end{split}
\end{align*}
where $\delta_{triv}$ denotes the connection homomorphism
corresponding to the sequence $0 \to A \to B \to C \to 0$ with the
trivial $W$-action. Thus we see $\delta=\delta_{triv} +
\mathbb{M}$ on cochains. Since general theory tells us that
$\delta$ and $\delta_{triv}$ give well defined homomorphisms
$H^n(W, C) \to H^n(W, A)$, we see that $\mathbb{M}$ does also. The
proof is complete.
\end{proof}

For computational purposes a more explicit form for $\mathbb{M}$
is desirable. To get this, note that the general cup product
construction gives us a map
$$H^1(W ,  \Hom(C,A)) \otimes H^n(W, C) \overset{\cup}{\to}
H^{n+1}(W, \Hom(C,A) \otimes C).$$ The composition pairing
$\Hom(C,A) \otimes C \to A$ induces a map $$H^{n+1}(W, \Hom(C,A)
\otimes C) \to H^{n+1}(W, A)$$ which when composed with the cup
product above yields a cup product:
$$H^1(W, \Hom(C, A)) \times H^n(W, C) \overset{\cup}{\to} H^{n+1}(W, A).$$

Finally note given a representation $\rho_B: W \to GL(B;A,C)$, we
have $$\log(\rho_B) \in \Hom(W, \Hom(C, A))=H^1(W, \Hom(C, A)).$$
Thus taking cup product with $\log(\rho_B) \in H^1(W, \Hom(C ,
A))$ induces a map
$$H^n(W, C) \to H^{n+1}(W, A).$$

Notice on cochains we have (see Brown \cite{Brown}, page 110)
$$(\log(\rho_B) \cup f)(z_0,\dots,z_n)
=(-1)^n \log ( \rho_B(z_0)) f(z_1,\dots,z_n),$$   and so we see that
the binding operator $\mathbb{M}$ is induced by cup product (from
the right) by $\log(\rho_B) \in H^1(W , \Hom(C, A))$. Thus it is
relatively routine to describe $\mathbb{M}$ in any computational
situation. We summarize this in the following proposition:

\begin{pro}
\label{pro:multwithlog} If $0 \to A \to B \to C \to 0$ is a short
exact sequence of $W$-modules with corresponding representation
$\rho_B: W \to GL(B; A,C)$ then considering $\log(\rho_B) \in
H^1(W, \Hom(C, A))$ we have $\mathbb{M}(-)= (-) \cup
\log(\rho_B)$.
\end{pro}

Now, we state our main result of this section.

\begin{thm}
\label{thm:mainresult} If $E(Q)$ lifts to an extension
with kernel $M$, then $\beta (q)= [\log ^* (\rho_M)]q$. \\
Conversely, if $\beta (q)=Lq$ for some $L \in H^1(W, \End(V))$,
then $E(Q)$ lifts to an extension with kernel $M$ where $M$ has
the representation $\exp^*(L) \in \Hom(W, GL(M; V,V))$.
\end{thm}

Note that in general $E(Q)$ may have more than one module where
the above lifting is possible. This is similar to having more than
one $L$ such that $\beta (q)=Lq$. Theorem~\ref{thm:mainresult}
says that choosing one fixes the other.

This result has many consequences for the structure of extension
classes. We will investigate them further in other sections. We
end this section with some calculations to illustrate the
effectiveness of the above result.

\begin{ex}
\label{ex:firstcalc} Let $W=V=(\FF _2)^ 3$ with standard basis and
$E(Q)$ be an extension with extension class $$q=(x_1^2 +x_2x_3,\
x_2 ^2+x_1 x_2,\ x_3 ^2+x_1 x_3).$$ Then, we have
$$\left[\begin{matrix}\beta(q_1)\cr \beta(q_2)\cr \beta(q_3)
\cr\end{matrix}\right]
=\left[\begin{matrix}0& x_3 & x_2\cr 0 & x_1 & 0 \cr 0 & 0 & x_1
\cr
\end{matrix}\right] \left[\begin{matrix}q_1\cr
q_2 \cr q_3\cr\end{matrix}\right].$$

This shows that $E(Q)$ has a lifting with kernel $M$ where $M$ has
representation $\rho : W \to GL(W;V,V)$ with
$$\rho (\alpha_1 , \alpha_2 ,  \alpha_3 )=
\left[\begin{matrix} 1 & 2 \alpha_3 & 2 \alpha _2\cr
0 & 1+ 2 \alpha _1 & 0 \cr 0 & 0 & 1+2 \alpha _1 \cr
\end{matrix}\right].$$
Note that $L$ is upper triangular, hence $M$ has a filtration
$$ 0 \leq \la (1,0,0) \ra \leq \la (1,0,0) , (0,1,0) \ra \leq M $$
whose sections are one dimensional $\ZZ/4 [W]$-lattices.
\end{ex}

\begin{ex}
\label{ex:secondcalc} Let $W=V= (\FF _2 )^3$ and $E(Q)$ be an
extension with extension class $$q=(x_1^2,\ x_2 ^2+x_1 x_2,\ x_3
^2+x_2 x_3).$$ Then, we have
$$\left[\begin{matrix}\beta(q_1)\cr \beta(q_2)\cr \beta(q_3)
\cr\end{matrix}\right]
=\left[\begin{matrix}0& 0 & 0 \cr 0 & x_2 & 0 \cr 0 & 0 & x_3 \cr
\end{matrix}\right] \left[\begin{matrix}q_1\cr
q_2 \cr q_3\cr\end{matrix}\right].$$ Thus, $E(Q)$ has a lifting
where in this case $M$ has the representation $$\rho (\alpha_1,
\alpha_2, \alpha_3 )= \left[\begin{matrix} 1 & 0 & 0 \cr 0 & 1+ 2
\alpha _2 & 0 \cr 0 & 0 & 1+2 \alpha _3 \cr
\end{matrix}\right].$$
In this case $L$ is diagonal, hence $M$ has a decomposition into one
dimensional lattices
$$M \iso  \la (1,0,0) \ra \oplus \la (0,1,0) \ra \oplus \la (0,0,1)
\ra .$$
\end{ex}

Note that applying Bockstein operator to the equation $\beta
(q)=Lq$, one gets that $[\beta(L)+L^2]q=0$. It turns out the
matrix $\beta (L) +L^2$ plays an important role. In the first
example above, we have
$$\beta(L)+L^2=
\left[\begin{matrix}0 & q_3 & q_2 \cr 0 & 0 & 0 \cr 0 & 0 & 0 \cr
\end{matrix}\right],$$
whereas we have $\beta (L)+L^2=0$ in the second case. We will see
in the next section that the matrix $\beta(L)+L^2$ is in fact the
obstruction for the module $M$ to decompose into a direct sum of
one dimensional lattices.

\section{Diagonalizable Extensions}
\label{sect:diagonalizability}

Let $E(Q): 0 \to V \to G(Q) \to W \to 0 $ be a Bockstein closed,
central extension corresponding to the quadratic map $Q: W \to V$.
Let $q\in H^2 (W, V)$ denote the extension class for $E(Q)$.

\begin{defn}
\label{defn:diagonalizable} We say that the quadratic map $Q$ is
diagonalizable if there exists a basis for $V$ such that the
components $q_1, \dots, q_n$ are individually Bockstein closed. In
other words, for some basis of $V$, we can write $q=(q_1, \dots,
q_n)$ where for each $i$ we have $\beta (q_i)=\lambda _i q_i$ for
some linear polynomial $\lambda _i$.
\end{defn}

It is clear that diagonalizable quadratic maps are Bockstein
closed with a diagonal binding matrix $L$. The converse also
holds:

\begin{lem}
\label{lem:conjugation} Let $Q: W \to V$ be a Bockstein closed
quadratic map with binding matrix $L$ for some basis of $V$. If
there is a invertible scalar matrix $N$ such that $N^{-1} L N$ is
a diagonal matrix, then $Q$ is diagonalizable.
\end{lem}

\begin{proof} Let $q'=N^{-1} q$. Then, $$\beta (q') =N^{-1} \beta
(q)= N^ {-1} L q= (N^{-1} L N )q'$$ shows that the coordinates of
$q'$ are individually Bockstein closed.
\end{proof}

When there is an invertible scalar matrix $N$ such that $N^{-1} L
N$ is a diagonal matrix, we say $L$ is diagonalizable. We have the
following useful criteria for the diagonalizability of $L$.

\begin{lem}
\label{lem:diagobstruction} Let $L$ be a matrix with linear
polynomial entries over the field $\FF _2$. Then, $L$ is
diagonalizable with a scalar matrix if and only if $\beta(L)+L^2=0$.
\end{lem}

\begin{proof} Write $L=\sum _i L_i x_i $. We have
$$\beta (L)+L^2=\sum _i (L_i^2+L_i) x_i^2 + \sum _{i<j}
(L_i L_j +L_j L_i) x_i x_j.$$ Thus, $\beta (L)+L^2 =0$ if and only
if  $L_i (L_i +I)=0$ for all $i$, and $L_i L_j = L_j L_i$ for all
$i,j$. Over the field $\FF _2$, a family $\{L_i\}$ is
simultaneously diagonalizable if and only if these equations are
satisfied. Note that if $N$ is the matrix that simultaneously
diagonalizes the family $\{ L_i \}$, then $N$ diagonalizes $L$ as
well.
\end{proof}

The argument used in the above proof is due to Dave Rusin. He uses
this argument to prove Lemma 20 in \cite{Rusin} which states that
all Bockstein closed extensions are diagonalizable. However, this
lemma is not correct. There are Bockstein closed extensions which
are not diagonalizable. The following is an example of such an
extension.
\begin{ex}
\label{ex:Rusincounterex} Let $E$ is an extension with $q=(x^2 +yz,\
y^2+xz,\ z^2)$. Then $$\beta(q)= (yz(y+z),\ xz(x+z),\ 0),$$ and we
can write
$$\left[\begin{matrix}\beta(q_1)\cr \beta(q_2)\cr \beta(q_3)\cr
\end{matrix}\right]
=\left[\begin{matrix}0& z & x+y\cr z & 0 & x+y \cr 0 & 0 & 0 \cr
\end{matrix}\right] \left[\begin{matrix}q_1\cr
q_2 \cr q_3\cr\end{matrix}\right].$$ So the equation $\beta (q)=Lq$
holds when $L$ is taken as the above coefficient matrix. In
particular, $E$ is a Bockstein closed extension. When we calculate
$\beta (L)+L^2 $, we get
$$\beta(L)+L^2=
\left[\begin{matrix}z^2 & z^2 & (x+y)(x+y+z)\cr z^2 & z^2 &
(x+y)(x+y+z) \cr 0 & 0 & 0 \cr
\end{matrix}\right] \neq 0,$$
so $L$ is not diagonalizable. This also implies that $E$ is not
diagonalizable since in this case the binding operator $L$ is
unique. Note that if $L'$ is another binding matrix satisfying
$\beta (q)=L'q$, then we would have $L=L'$ since the components
$$x^2 +yz,\ y^2+xz,\ z^2$$ of $q$ form a regular sequence. We
conclude that $E$ is a Bockstein closed extension which is not
diagonalizable.
\end{ex}

Now, we continue to find more equivalent conditions for
the diagonalizability of $L$.
\begin{pro}
\label{pro:liftinglattices} Let $M$ be a $\ZZ/4 [W]$-lattice with
representation $\rho _M \in GL (M; V, V)$, and let $L=\log \rho _M
\in H^1 (W, \End (V))$. Then,
the following are equivalent: \\
(i) $L$ is diagonalizable, i.e., $\beta(L)+L^2=0$. \\
(ii) $M$ is a direct sum of one-dimensional $\ZZ /4[W]$-lattices.\\
(iii) $M$ lifts to a $\ZZ [W]$-lattice. \\
(iv) $M$ lifts to a $\ZZ /8[W]$-lattice.
\end{pro}

\begin{proof}
$(i) \Leftrightarrow (ii)$ follows from the fact that the
representation $\rho: W \to GL(M)$ of $M$ is defined as $\rho
(w)=1+2L(w)$ mod $2$. To see that $(ii)\Rightarrow (iii)$, note that
we just need to show that every one dimensional $\ZZ / 4[W]$-lattice
lifts to an integral lattice. This follows from the fact that the
unit group of $\ZZ /4$ is $\{\pm 1\}$ which is also the unit group
of $\ZZ$. The implication $(iii) \Rightarrow (iv)$ is obvious. To
complete the proof, we will show that $(iv) \Rightarrow (i)$.

Let $\tilde \rho: W \to GL(\widetilde M)$ denote the lifting of
$\rho$. We can write $\tilde \rho (w)=I+2L(w)+4a(w)$ mod $8$ for
some $a : W \to \End(V)$. Since $2w=0$, we have
$$ I=[\tilde \rho (w)]^2=I + 4 [L(w)^2+ L(w)] \ \ {\rm mod} \ \ 8 $$
hence $L(w)^2+L(w)=0$ mod 2 for every $w \in W$. This gives
$BL+L^2=0$.
\end{proof}

As a consequence we obtain the following interesting result:

\begin{pro}
\label{pro:mainthmforlifting} Let $E': 0 \to M \to \Gamma \to W
\to 0$ be an extension with $M$ a $\ZZ$-free $\ZZ [W]$-module such
that $M /2M$ is a trivial $\ZZ /2 [W]$-module. Then, the extension
$$ E: 0 \to M/2M \to \Gamma / 2M \to W \to 0$$
is a diagonalizable extension.
\end{pro}

\begin{proof} Let $Q$ be the quadratic map for the
extension $E$. Since $E$ has a lifting
$$ E'' : 0 \to M/ 4M \to \Gamma / 4M \to W \to 0$$
it is a Bockstein closed extension whose binding operator $L$ is
determined by $M /4M$. Since the $\ZZ /4[W]$-lattice $M/4M$ is
reduced from an integral module, by
Proposition~\ref{pro:liftinglattices}, $L$ is diagonalizable, so $Q$
is diagonalizable as well.
\end{proof}

The converse of Proposition~\ref{pro:mainthmforlifting} is also
known to be true:

\begin{pro}
\label{pro:integrallifting} Let $E : 0 \to V \to G \to W \to 0$ be
a diagonalizable extension. Then, it lifts to a (unique) extension
$ \widetilde E : 0 \to M \to \widetilde G \to W \to 0$ where $M$
is a $\ZZ $-free $\ZZ [W]$-module.
\end{pro}

To prove this we will need the following lemma:

\begin{lem}
\label{lem:cutelemma} Let $q$ be a quadratic polynomial in $m$
variables viewed as an element in $H^2(W,\FF _2)$ where
$m=\dim(W)$ as before. Then $q$ is reducible if and only if $\beta
(q)= \ell q$ for some linear polynomial $\ell \in H^1(W,\FF _2)$.
In this case if $q=uv$, then in fact $\ell =u+v$.
\end{lem}
\begin{proof}
If $q$ is reducible with $q=uv$ a simple calculation shows $\beta
(q)=\ell q$ with $\ell =u+v$.

Let us prove the converse so assume $\beta (q)= \ell q$ for some
linear polynomial $\ell$. We will show that this implies that $q$ is
decomposable. If $m=1$ there is nothing to show, so assume $m > 1$.

Case 1: $\ell =0$: In this case, $\beta (q)=0$ and so $q=\beta
(u)=u^2$ by the fact that $H^*(W, \FF _2)$ is $\beta$-acyclic.

Case 2: $\ell \neq 0$: In this case as $\ell $ is nonzero in
$H^1(W,\FF _2)=\Hom(W,\FF _2)$ we can let $H=\ker( \ell )$ and $H$
will be a hyperplane in $W$. Thus $\beta (q)=0$ (and hence $q$ is
a square) when restricted to $H^*(H, \FF _2)$ and so it follows
that $q=u^2 + \ell v \in H^*(W,\FF _2)$ where $u,v$ are linear and
$u$ is algebraically independent from $\ell$. Applying $\beta$ to
both sides of the last equation we get $\ell q=\ell ^2v + \ell
v^2$. Canceling $\ell$ we get $q=\ell v+v^2=v(\ell +v)$ and so $q$
is decomposable as desired.
\end{proof}

We will now prove Proposition~\ref{pro:integrallifting}:
\begin{proof}
Take a basis for $V$ such that, as a $\ZZ/2 [W]$-module, $V$
decomposes as $V = \oplus_{i=1}^n V_i$ where $\dim(V_i)=1$. Taking
components of the extension class $q$ with respect to this basis,
say $\{q_1,\dots,q_n\}$, we find that $q_i$ represents an
extension
$$0 \to V_i \to G_i \to W \to 0$$
and the $q_i$ are individually diagonal, i.e., $\beta (q_i)=\ell
_iq_i$.

Suppose that we can uniformly lift the corresponding extensions
for the individual $q_i$, say to extensions
$$ 0 \to M_i \to \widetilde{G}_i \to W \to 0$$
with extension class $\tilde{q}_i$, where $M_i$'s are $\ZZ
/4[W]$-modules with $M_i/2M_i=V_i$ and $M_i \cong \ZZ /4$ as
abelian groups. The isomorphism $$ H^2 (W,\oplus_{i=1}^n M_i)
\cong \oplus_{i=1}^n H^2(W,M_i)$$ shows that the $\hat{q}_i$ would
fit together to give a uniform lift $\tilde{q}=(\tilde{q}_1,
\dots, \tilde{q}_n)$ with lifting module $M=\oplus_{i=1}^n M_i$.
Thus it is sufficient to prove the theorem in the case
$\dim(V)=1$.

Since $\beta (q)=\ell q$ by assumption, Lemma~\ref{lem:cutelemma}
gives that $q=uv$ for some 1-dimensional classes $u,v$. There are
then two cases: \\
Case 1: $\{u,v\}$ are linearly dependent. In this case $G$ is abelian
and a uniform integral lift certainly exists. \\
Case 2: $\{u,v\}$ are linearly independent. In this case it is
easy to check that $G=D_8 \times S$ where $S$ is elementary
abelian with $\dim(S)=\dim(W)-2$ and $V=[D_8,D_8]$. This extension
has an integral lift (unique) given by $\widetilde{G}=D_{\infty}
\times S$ and $M=[D_{\infty},D_{\infty}]$ where $D_{\infty}$ is
the infinite dihedral group given by a semidirect product of $\ZZ$
by $\ZZ /2$ with twisting given by the sign map.
\end{proof}

We end this section with a result which summarizes the results
obtained about diagonalizability.

\begin{thm}
\label{thm:summary} Let $Q : W \to V$ be a Bockstein closed
quadratic map, and let $E: 0 \to V \to G \to W \to 0 $ be the
central extension associated to $Q$. Let $q$ denote the extension
class for $E$. Then, the following are equivalent: \\
(i) $Q$ is diagonalizable, i.e., there is a basis for $V$ such that
the components of $q$ are individually
Bockstein closed. \\
(ii) There is a choice of basis of $V$, such that the components
$q_1, \dots, q_n$ of $q$ all decompose as $q_i=u_iv_i$
where $u_i,v_i$ are linear polynomials. \\
(iii) There exists a diagonalizable $L \in \Hom (W, \End (V))$ such
that $\beta(q)=Lq$.
This is characterized exactly by the equation $\beta(L)+L^2=0$. \\
(iv) $E$ lifts to an extension with kernel $M$ where $M$ is a
direct sum of one dimensional $\ZZ /4 [W]$-lattices. \\
(v) $E$ lifts to an extension with kernel $M$ where
$M$ is a $\ZZ /8 [W]$-lattice. \\
(vi) $E$ lifts to an extension with kernel $M$ where $M$ is a $\ZZ
[W]$-lattice, i.e., $E$ has a uniform integral lifting.
\end{thm}

\begin{proof} $(i)\Leftrightarrow (ii)$ follows from Lemma~\ref{lem:cutelemma},
and Lemma~\ref{lem:conjugation} gives $(i)\siff (iii)$. We also
have $(i) \siff (vi)$ by Propositions~\ref{pro:mainthmforlifting}
and~\ref{pro:integrallifting}. Finally, the equivalences $(iv)
\siff (v) \siff (vi)$ follow from
Proposition~\ref{pro:liftinglattices}.
\end{proof}

\section{Triangulable Extensions}
\label{sect:triangulable}

As in the previous section, $E(Q): 0 \to V \to G(Q) \to W \to 0$
denotes an arbitrary Bockstein closed extension associated to a
quadratic map $Q: W \to V$. Let $q \in H^2 (W, V)$ be the
associated extension class.

\begin{defn}
\label{def:triangulable} We say the quadratic map $Q$ is (upper)
triangulable if there is a basis for $V$ such that the components
$q_1, \dots, q_n $ of $q$ have the property that for each
$i=1,\dots, n$, the ideal $(q_i, q_{i+1}, \dots , q_n )$ is a
Bockstein closed ideal.
\end{defn}

Note that if $Q$ is triangulable then it is Bockstein closed with
an upper triangular binding matrix $L$. The converse also holds:
If $Q$ is a quadratic map with a binding matrix $L$ for some
basis, then $Q$ is triangulable if there is a scalar matrix $N$
such that $N^{-1} L N $ is a upper triangular matrix. In this case
we say $L$ is triangulable.

>From our earlier discussion about the connection between $L$ and
the lifting lattice $M$, the following is immediate.

\begin{lem}
\label{lem:triangulable} Let $M$ be a $\ZZ/4 [W]$-lattice with
representation $\rho _M \in GL (M; V, V)$, and let $L=\log \rho _M
\in H^1 (W, \End (V))$. Then, $L$ is triangulable if and only if $M$
has a filtration $0 \subseteq M_1 \subseteq M_2 \subseteq \dots
\subseteq M_n=M$ such that each factor $M_i/M_{i-1}$ is a rank one
$\ZZ /4[W]$-lattice.
\end{lem}

There are many $\ZZ /4[W]$-lattices $M$ which do not have such a
filtration even with the extra condition that $M/2M$ is a trivial
$\ZZ /2 [W]$-module. Note that given a $\ZZ /4[W]$-lattice $M$ such
that $V=M /2M$ is trivial, there is a homomorphism  $$L : W\to \End
(V)$$ associated to it under the exponential-logarithm
correspondence. It is easy to see that if $M$ has a filtration with
one dimensional factors, then the family $\{ L(w) \ | \ w \in W \}$
is simultaneously triangulable. This means that there is a matrix
$N$ such that $N^{-1} L (w) N =T(w)$ is an upper triangular matrix
for all $w \in W$. Since $T(w)^2 +T(w)$ is a strictly upper
triangular matrix, it follows that $[T(w)^2 +T(w)]^n=0$ where
$n=\dim V$. This implies that $[L(w)^2 +L(w)]^n=0$ as well. Note
that for $w=\alpha_1 w_1 + \cdots + \alpha _m w_m $, the scalar
matrix $L(w)^2 +L(w)$ is equal to the value of $\beta (L )+L^2 $
calculated by setting $x_1= \alpha _1,\ x_2=\alpha _2, \dots,\ x_n
=\alpha _n $. So, if $[L(w)^2 +L(w)]^n=0 $ for all $w \in W$, then
we have $[\beta (L) +L^2 ]^n =0$. We have proved the following:

\begin{lem}
\label{lem:trigcondition} If $L$ is triangulable, then $\beta (L )
+L^2 $ is nilpotent.
\end{lem}

Using this we can give an example of a $\ZZ /4[W]$-lattice which
does not have a filtration:

\begin{ex}
\label{ex:nofiltration} Let $m=n=2$ and $\{ w_1, w_2\}$ be a basis
for $W$. Consider the representation $\rho : W \to GL_2 (\ZZ /4)$
where
$$\rho (w_1)=\left[\begin{matrix} 1 & 2 \cr 0 & 1 \cr
\end{matrix}\right], \ \ \ \ \ \  \rho(w_2)=\left[\begin{matrix} 3 & 0
\cr
2 & 1 \cr \end{matrix}\right].$$
Then, we have
$$L=\left[\begin{matrix} x_2 & x_1 \cr x_2 & 0  \cr
\end{matrix}\right] \ \  {\rm and}  \ \ \ \  \beta (L)+L^2=
\left[\begin{matrix} x_1x_2 & x_1^2+x_1x_2 \cr 0 & x_1 x_2 \cr
\end{matrix}\right].$$ It is clear that $[\beta (L) +L^2]^k \neq
0$ for any $k$, so $M$ does not have a one dimensional sublattice.
\end{ex}

In the previous section, we showed that $L$ is diagonalizable if
and only if $\beta (L)+L^2 =0$. So, it is reasonable to ask if the
converse of Lemma~\ref{lem:trigcondition} holds. As positive
evidence one sees that if $\beta (L)+L^2$ is nilpotent then for
every $w\in W$, the operator $L(w)$ is triangulable. This is
because, if $\beta (L) +L^2 $ is nilpotent, i.e., $[\beta (L)
+L]^k =0$ for some $k$, then
$$[L(w)^2 +L(w)]^k= L(w)^k [L(w)+I]^k =0$$ holds for every $w \in W$.
So, the minimal polynomial of $L(w)$ is a product of linear
polynomials, and hence $L(w)$ is triangulable for all $w \in W$ by a
standard result in linear algebra. Unfortunately, unless we have an
extra structure, in general we do not have simultaneous
triangulability. In fact, the following example clearly shows that
the converse of Lemma~\ref{lem:trigcondition} fails.

\begin{ex}
\label{ex:devil} Let $m=n=2$ and $\{ w_1, w_2\}$ be a basis for
$W$. Consider the representation $\rho : W \to \Aut (M)=GL_2 (\ZZ
/4)$ where
$$\rho (w_1)=\left[\begin{matrix} 1 & 2 \cr 0 & 1 \cr
\end{matrix}\right], \ \ \ \ \ \  \rho(w_2)=\left[\begin{matrix} 1 & 0
\cr
2 & 1 \cr \end{matrix}\right].$$ Then,
$$L=\left[\begin{matrix} 0 & x_1 \cr x_2 & 0  \cr
\end{matrix}\right] \ \  {\rm and}  \ \ \ \  \beta (L)+L^2=
\left[\begin{matrix} x_1x_2 & x_1^2  \cr x_2 ^2 & x_1 x_2 \cr
\end{matrix}\right].$$ It is clear that $[\beta (L) +L^2]^2=0$. It
is easy to check that no nonzero vector is a common $\ZZ
/4$-eigenvector for $\rho(w_1)$ and $\rho(w_2)$ and thus $M$ has no
one dimensional $\ZZ /4 [W]$-sublattice.
\end{ex}

The above examples show that the situation with triangulability is
much more complicated. To illustrate that such bad examples also
appear as extensions, we calculate $\beta (L) +L^2 $ for the
$\gl$-induced extension given in Example~\ref{ex:gl}.

\begin{ex}
\label{ex:glagain} Consider the central extension $$E(Q _{\gl
_n}): 0 \to \gl _n (\FF _2 ) \to G(Q_{\gl _n}) \to \gl _n (\FF _2
) \to 0
$$ with quadratic map $Q _{\gl _n} (\MA) = \MA + \MA ^2$. Note that
we can express the extension
class
$$q \in H^2 ( \gl _n (\FF _2 ), \gl _n (\FF _2 ) )$$ as a $n
\times n $-matrix $\MQ$ whose $ij$-th entry will be the $ij$-th
component of $q$. Then one computes
$$ \MQ _{ij} (\MA)= (\MA ^2 + \MA )_{ij} = \sum _k a_{ik} a_{kj} + a_{ij}$$
So, we can write $q_{ij} = x_{ij} ^2 + \sum _k x_{ik} x_{kj}$
where $x_{ij} \in H^1 ( \gl _n (\FF _2 ) ,  \ZZ /2 )$ is the
linear form which takes a matrix $\MA$ to its $ij$-th entry. Note
that if we set $$\MX \in H^1 ( \gl _n (\FF _2 ) , \gl _n (\FF _2 )
)$$ as the class corresponding to the identity homomorphism $id:
\gl _n (\FF _2 ) \to \gl _n (\FF _2 )$, then $\MX$ will be a
matrix with $ij$-th entry $x_{ij}$, and we will have
$$\MQ =\beta (\MX) + \MX^2.$$ Here the ring structure on
$H^*(\gl _n (\FF _2 ) ,\gl _n (\FF _2 ) )$ that we are using is
induced from the composition map $\gl _n (\FF _2 ) \times \gl _n
(\FF _ 2) \to \gl _n (\FF _ 2) $ and is noncommutative!

This gives us an easy way to calculate the Bockstein from which we
obtain
$$ \beta ( \MQ )= \beta (\MX) \MX + \MX \beta (\MX)= \MQ \MX+ \MX \MQ.$$
This shows that for all $(i,j)$ pairs, $\beta (q_{ij})$ lies in
the ideal generated by the components of $q$, and hence provides
another way to see that the quadratic map $Q_{\gl _n}$ is
Bockstein closed. Note that the matrix $L$ with respect to some
basis can be written as an $n^2\times n^2$-matrix, but it is much
more convenient to think of $L$ as a homomorphism  $$L : \gl _n
(\FF _2 ) \to \End ( \gl _n (\FF _2 ) ))$$ such that for every
$\MB \in \gl _n (\FF _2 ) $ the image of $\MB$ under $L$ is
defined as the endomorphism
$$L(\MB ) : \MA \to [\MA, \MB]=\MA \MB +\MB \MA.
$$ Using $L(\MX ): \MA \to [\MA, \MX]$, we can calculate
$$[\beta (L (\MX))
+L(\MX) ^2]: \MA \to [\MA,\ \beta (\MX) +\MX^2]=[\MA, \MQ].$$ This
shows, in particular, that $$[\beta (L)+L^2 ]q=[\beta(L(\MX) )
+L(\MX)^2 ] (\MQ)= [\MQ,\MQ]=0$$ which we know holds for all
extension classes. Note that this extension is not triangulable
because $\MQ$ is not a nilpotent matrix.

On the other hand, if we had taken $\mathfrak{u} _n$ instead of $\gl
_n$, then we would have $\MQ= \beta (\MU) +\MU ^2$ where $\MU$ is a
strictly upper triangular matrix with $ij$-th entry equal to
$x_{ij}$. It is clear that $\MQ$ is also strictly upper triangular,
so $\MQ ^n$ will be zero. In fact, $Q_{\mathfrak{u} _n}$ is a
triangulable quadratic map.
\end{ex}

\begin{rem}
\label{rem:liftingobstruction} The calculation we performed above
in Example \ref{ex:glagain} can be used to see some of the earlier
results in a different way. Let $E: 0 \to V \to G \to W \to 0$ be
a Bockstein closed extension with binding operator $L$. Choosing a
basis for $V$, we can view $L$ as a homomorphism $L : W \to \gl _n
(\FF _2 )$. Using $L$, we can lift the extension $E(Q_{\gl _n} )$
and obtain
$$
\begin{CD}
E(Q_L):\ @. 0 @>>> \gl _n (\FF _2 ) @>>> G (Q_L) @>>> W @>>> 0\\
@. @. @| @VVV @VVLV \\
E(Q_{\gl _n}):\ @. 0 @>>> \gl _n (\FF _2 ) @>>> G (Q_{\gl _n })
@>{\pi}>> \gl _n (\FF _2 ) @>>> 0.\\
\end{CD}
$$
Since the extension class for the bottom extension is $\beta
(\MX)+\MX ^2$, the extension class for $E(Q_L)$ will be $\beta
(L)+L^2$. So, the class $\beta (L)+L^2$ can be thought of as the
obstruction for lifting $L: W \to \gl _n (\FF _2 )$ to a
homomorphism $\widehat L : W \to G (Q_{\gl _n })$ such that $\pi
\circ \widehat L  =L$. Note that $L$ corresponds to a group
homomorphism $\rho : W \to K_n (\ZZ /4)$ under the exp-log
correspondence, and the map $\pi $ becomes the mod $4$ reduction map
$$K_n (\ZZ /8 ) \to K_n (\ZZ /4).$$
when $K_n (\ZZ /4)$ is identified with $\gl _n (\FF _2 )$ (see
Example \ref{ex:gl}). This shows that $\beta (L) +L^2$ is the
obstruction for lifting the representation $\rho : W \to K_n (\ZZ
/4)$ to a representation of a $\ZZ /8[W]$-lattice. This provides
another way to see the equivalence $(i)\Leftrightarrow (iv)$ given
in Proposition \ref{pro:liftinglattices}.
\end{rem}

\section{Frattini and Effective Extensions}
\label{sect:Frattini&Effective}

In this section we consider Bockstein closed extensions with some
additional conditions. These conditions are standard conditions one
considers in group extension theory. We start with the following:

\begin{defn}
\label{def:Frattiniext} A central extension $E: 0 \to V \to G \to
W \to 0 $ is called a {\bf Frattini extension} if the Frattini
subgroup $\Phi(G)$ of $G$ is equal to $V$.
\end{defn}

Recall that the  Frattini subgroup of a $p$-group $G$ is defined
as the subgroup $\Phi(G)=G^p [G,G]$ where $G^p$ is the subgroup
generated by the $p$th powers and $[G,G]$ is the commutator
subgroup of $G$. In the case $p=2$, one has $\Phi (G) =G^2$, since
in this case, $G/G^2$ has exponent $2$ implies that it is abelian
so $[G,G] \leq G^2$. Because of this it makes sense to define the
Frattini condition for quadratic maps as follows.

\begin{defn}
\label{def:FrattiniQ} We say a quadratic map $Q : W \to V$ is {\bf
Frattini} if the set $\Im(Q)=\{ Q(w) \ | \ w \in W \}$ generates
$V$. Equivalently, a quadratic map $Q$ is Frattini if the
associated central extension $E(Q) : 0 \to V \to G(Q) \to W \to 0
$ is a Frattini extension.
\end{defn}

Notice that when an extension of the form $E: 0 \to V \to G \to W
\to 0$ is not Frattini, then $\Phi(G)$ is a proper subspace of
$V$, and the group $G$ splits as $G' \times \ZZ /2$. In terms of
group extension theory the trivial summand causes no extra
difficulties, so to avoid trivialities, one often assumes that the
extensions in question are Frattini extensions.

For Frattini extensions we have the following useful criteria:

\begin{lem}
\label{lem:Frattini} Let $E: 0 \to V \to G \to W \to 0$ be a
Frattini extension with extension class $q=(q_1, \dots, q_n)$ with
respect to some basis for $V$. If \ $ a_1q_1+\cdots + a_n q_n =0$
for some scalars $a_1,\dots, a_n$, then $a_1=\dots =a_n=0$.
\end{lem}

\begin{proof} If there exists $a_1, \dots, a_n \in \FF _2$ with
$a_1q_1+\cdots + a_n q_n =0$, then one would have $a_1Q_1(w)+\cdots
+ a_n Q_n (w) =0$ for all $w \in W$. But, then $\Im (Q)$ will lie inside
the kernel of the functional $a_1 x_1 + \dots + a_n x_n $. The Frattini
condition $\Im (Q)=V$ will hold only if this functional is zero.
\end{proof}

We also have the following nice basis choice for Frattini
extensions.

\begin{lem}
\label{lem:nicebasis} Let $Q : W \to V$ be a Frattini extension
where $W$ and $V$ are nonzero, and let $k$ be a positive integer
such that $k \leq \min \{ \dim V , \dim W \}$. Then, there exists
a set of linearly independent vectors $\{ w_1,\dots, w_k \}$ in
$W$ such that $\{ Q(w_1),\dots , Q(w_k)\}$ is a linearly
independent set in $V$.
\end{lem}

\begin{proof}
Let $m=\min\{ \dim V, \dim W \}$. For $k=1$, the lemma is clear,
since there is always a non-zero vector $w$ in $W$ whose image is
non-zero, otherwise $V$ will be the zero space. Assume that the
lemma is true for some $k<m$. We will show that it also holds for
$k+1$.

By assumption there is a linear independent set $\{ w_1, \dots,
w_k\}$ of vectors in $W$  such that $\{ Q(w_1), \dots, Q(w_k)\}$
is a linearly independent set in $V$.  Let $W_k$ be the subspace
of $W$ generated by $\{ w_1,\dots, w_k \}$, and $V_k$ the subspace
of $V$ generated by $\{ Q(w_1),\dots , Q(w_k)\}$. Then, $\dim W_k
=\dim V_k =k <m$. To complete the proof we need to show that there
exists a $w$ in the set difference $W-W_k$ such that $Q(w) \not
\in V_k$. Assume to the contrary that there is no such $w \in W$,
i.e., $Q(w)$ lies in $V_k$ for all $w \in W-W_k$. Now fix a $w \in
W-W_k$. Note that for all $i=1,\dots , k$, we have $Q(w+w_i)\in
V_k$. Therefore, for all $1\leq i \leq k$, we have
$$B(w, w_i)=Q(w+w_i)+Q(w)+Q(w_i) \in V_k.$$
By bilinearity, we get $B(w, w')\in V_k$ for all $w' \in W_k$.
Since this statement is true for all $w \in W-W_k$, we have $$B(
W-W_k, W_k) \subseteq V_k.$$

Now, take $w \in W-W_k$ and $w' \in W_k$. Since $B(w,w')$, $Q(w)$,
and $Q(w+w')$ are all in $V_k$, we can conclude that
$$ Q(w')=B(w, w')+Q(w)+Q(w+w')$$
is in $V_k$. Thus $Q(w') \in V_k$ for every $w' \in W_k$. For
vectors in $W-W_k$ we assumed at the beginning that their image
under $Q$ lies in $V_k $, so we obtained that $Q(w)$ lies in $V_k$
for all $w \in W$. But $V_k$ is a proper subspace of $V$ since
$\dim V_k=k <m\leq \dim V$. This contradicts the assumption that
$Q$ is Frattini.
\end{proof}

This last lemma, in particular, tells us that if $\dim V=\dim
W=n$, then there exists a basis $\{w_1, \dots, w_n \}$ for $W$ and
a basis $\{ v_1, \dots, v_n \}$ for $V$ such that $Q(w_i)=v_i $
for all $i=1, \dots, n$. In this case, the extensions class
$q=(q_1, \dots, q_n)$ is in the form
$$ q_k =x_k ^2 + \sum _{i<j} \gamma ^{(k)} _{ij} x_i x_j $$
for all $k=1, \dots, n$. If the extension class is written in this
form we say it is in {\bf bijective form} and the basis which is
induced from this form is called a {\bf bijective basis}.

Note that in the case where $\dim W=\dim V=n$, a bijective basis
also allows us to identify $W$ with $V$ and write quadratic maps
as operators on one vector space. We will use this later in the
paper.

Now, we will impose another property for the extensions $E(Q)$.

\begin{defn}
\label{defn:effective} We say a quadratic map $Q : W \to V$ is
{\bf effective} if $$\Ker (Q)=\{ w \in W  \ |\ Q(w)=0 \}=\{ 0\}.$$
An extension $E(Q): 0 \to V \to G(Q) \to W \to 0$ is called
effective if the associated quadratic map $Q$ is effective.
\end{defn}

Equivalent interpretations of effective extensions are given in
the following proposition:

\begin{pro}
\label{pro:effectiveness} Let $E(Q): 1 \to V \to G(Q) \to W \to 1$
be a central extension with quadratic map $Q$ and extension class
$q$. Then the following are equivalent:

(i) $E(Q)$ is effective, i.e., $Q(w)=0 $ implies $w=0$ for all $w
\in W$.

(ii) $\Res ^W_{\la w \ra } (q) \neq 0$ for any nonzero cyclic
subgroup $\la w \ra$ of $W$.

(iii) $V$ is a maximal elementary abelian subgroup of $G(Q)$.

(iv) All the elements of order $2$ in $G(Q)$ lie in $V$.
\end{pro}

\begin{proof} It is easy to see that $(i) \Leftrightarrow (ii)$ by using
the fact that any factor set $f$ representing $q$ has $f(w,w)=Q(w)$.
Note that on a one dimensional subspace $\la w \ra$ of $W$, the
vector $Q(w)=f(w,w)$ determines whether $Q$ and $q$ are zero when
restricted to that subspace. The implication $(i) \Rightarrow (iii)$
is clear since if $V$ is not maximal then one can find non trivial
element $w \in W$ coming from the larger elementary abelian subgroup
such that $Q(w) = (\hat w ) ^2 =0$. On the other hand if we have
$(iii)$ and $g \in G$ has order $2$, since $V$ is central in $G$, we
have $ \la g ,V \ra$ is an elementary abelian $2$-group. Since $V$
is maximal, it must be that $g \in V$. Thus $(iii) \Rightarrow(iv)$.
Finally, $(iv) \Rightarrow (i)$ since if there is a lift $\hat{w}
\in G-V$ such that $(\hat{w})^2=1$, then there exists a nonzero $w
\in W$ such that $Q(w)=0$.
\end{proof}

Note since a nontrivial zero of $Q: W \to V$ corresponds to a
nontrivial common zero of the homogeneous quadratic polynomials
$Q_1,\dots,Q_m: W \to \FF _2 $, the Chevalley-Warning theorem tells
us that if $E$ is an effective extension then $\dim W \leq 2 \dim
V$. This bound can be seen to be strict by considering the effective
central extension for $G=(Q_8)^m$ where $Q_8$ is the quaternionic
group of order 8 given by
$$0 \to (\ZZ /2)^m \to (Q_8)^m \to (\ZZ /2)^{2m} \to 0$$
where $\dim W=2m$ and $ \dim V=m$.

For $Q: \gl _n (\FF _2)  \to \gl _n (\FF _2 ) $, one has
$Q(\MA)=\mathbb{O}$ if and only if $\MA ^2=\MA$ (i.e., $\MA$ is a
projection). Thus if $W$ is a square-closed subspace, $Q: W \to W$
is effective whenever $W$ does not contain any nonzero
projections. Thus for example $W=\mathfrak{u} _n$, the strictly
upper triangular matrices. (Since the only projection with all
zero eigenvalues is the zero projection.)

\begin{defn}
\label{defn:regularsequence} Let $f_1, \dots, f_m$ be a sequence
of polynomials in $\FF _2 [x_1,\dots, x_n]$. The sequence $f_1,
\dots , f_m$ is called a regular sequence if $f_k$ is a non-zero
divisor in the quotient ring $\FF _2 [x_1, \dots , x_n]/(f_1,
\dots , f_{k-1})$ for all $k=1,\dots,n$.
\end{defn}

We have the following important fact about extensions which are
both Bockstein closed and effective:

\begin{pro}[Serre, Quillen, Carlsson]
\label{pro:regularsequence} Let $E$ be a Bockstein closed, effective
extension. Then, $\dim W \leq \dim V$. Moreover, if $\dim V=\dim W$,
then the components $q_1, \dots, q_n$ of the extension class  (and
any permutation of them) form a regular sequence.
\end{pro}

\begin{proof} Let $k$ denote the algebraic closure of $\FF _2 $.
When $E$ is Bockstein closed, the associated ideal $I(Q)$
generated by the components of the quadratic map is closed under
the Steenrod operations. Then, by a result of Serre \cite{Serre},
the variety of $I(Q)$ over $k$ must include a non-trivial rational
point if it is a nonzero variety. (Note since $I(Q)$ is a
homogeneous ideal, zero will always be in its corresponding
variety. Serre's theorem says that if the variety has some other
nontrivial solution over $k$, it will also have a nontrivial
solution over $\FF _2$.) By the effectiveness condition, the
components $q_1,\dots, q_n$ of the extension class $q$ have no
nontrivial common zero in $W$. So, $q_1,\dots, q_n$ have no
nontrivial common zero in $k$, by Serre's theorem. So, $n \leq m$
where $m=\dim W$ and $n=\dim V$ as before.

For the second part we use a result in commutative algebra which
states that a sequence of polynomials $f_1,\dots, f_n$ in $m$
variables is a regular sequence if
$$\dim Var(f_1,\dots, f_m) \leq m-n$$ (see for example \cite{Quillen}).
Since $m-n=0 =\dim {\rm Var}(q_1,\dots , q_n)$, we can conclude
that $q_1,\dots, q_n$ (or any permutation of them) is a regular
sequence.
\end{proof}

Thus we are motivated to make the following definition:

\begin{defn}
\label{defn:2-powerexact} A central extension
$$E(Q): 0 \to V \to G(Q) \to W \to 0$$
with corresponding quadratic map $Q: W \to V$ is called
2-power exact if the following 3 conditions hold: \\
(a) $\dim(V)=\dim(W)$. \\
(b) The extension is Frattini, i.e., $\Im(Q)=V$. \\
(c) The extension is effective, i.e., $Q(w)=0$ if and only if
$w=0$.
\end{defn}

In the next section we will study Bockstein closed $2$-power exact
extensions. Note that if $E$ is a Bockstein closed extension
satisfying only the conditions $(a)$ and $(c)$ above then we can
conclude that it is $2$-power exact. To see this, observe that the
quadratic map $Q: W \to \Im (Q)$ is Bockstein closed and
effective, so by Proposition \ref{pro:regularsequence}, we should
have $\dim W \leq \dim \Im (Q)$. By condition $(a)$, we have $\dim
W=\dim V$, so we can conclude $\Im (Q)=V$.

Note that a typical example of a $2$-power exact extension is the
extension associated with the $\gl$-induced quadratic map $Q:
\mathfrak{u} _n (\FF _2 ) \to \mathfrak{u} _n (\FF _2 ) $ where
$\mathfrak{u} _n (\FF _2)$ is the $\FF _2$-vector space of
strictly upper triangular matrices.

\section{Bockstein Closed $2$-Power Exact Sequences}
\label{sect:2-powerexact}

In this section we consider Bockstein closed $2$-power exact
extensions and prove some restrictions on binding operators of these
extensions. This leads to an interesting group theoretical result
for such extensions. We start with the following:

\begin{pro}
\label{pro:uniqueness} Let $E(Q): 0 \to V \to G(Q) \to W \to 0$ be a
Bockstein closed $2$-power exact sequence corresponding to the
quadratic map $Q: W \to V$. Then, there is a unique binding operator
for $Q$.
\end{pro}

\begin{proof} Let $L_1, L_2 \in H^1(W, \End(V))$ such that
$\beta (q)=L_1 q =L_2q$. Then we have $(L_1+L_2)q=0$. Choosing a
basis for $W$ and $V$, we can express $q$ as a column vector and
$L$ as a $n\times n$ matrix with entries in $H^1(W, \FF _2)$. So,
equation $(L_1 +L_2 )q=0$ gives a system of equations in $q_i$'s
with coefficients in $H^1(W, \FF _2)$. By
Proposition~\ref{pro:regularsequence}, the entries of $L_1+L_2$
must lie in $I(Q)$. Since the entries are one dimensional and
$I(Q)$ is generated by $2$-dimensional classes, we have
$L_1+L_2=0$.
\end{proof}

Another consequence of regularity is the following:

\begin{pro}
\label{pro:specialform} Let $E(Q): 0 \to V \to G(Q) \to W \to 0$
be a Bockstein closed 2-power exact sequence corresponding to the
quadratic map $Q: W \to V$. Choose a bijective basis for $V$ and
$W$ as in Lemma~\ref{lem:nicebasis}, and let  $ L = \sum _i L_i
x_i$ be the binding matrix with respect to this basis. Then, we
have
$$ \beta (L) + L^2 = \sum _i (L_i ^2 + L_i ) q_i.$$
\end{pro}

\begin{proof} Applying the  Bockstein operator on $\beta (q)=Lq$,
we get $[\beta (L ) +L ^2 ]q=0$. Since $q_1, \dots, q_n $ is a
regular sequence (in any order), each entry of $\beta (L ) +L^2 $
must be in the ideal $I(q)$. So, we have $$ \beta (L) + L^2 = \sum
_i K_i q_i $$ for some $n \times n$-matrices $K_1, K_2, \dots, K_n
$. Note that $$\beta (L) +L ^2 =\sum _i (L_i ^2 + L_i ) x_i ^2 +
\sum _{i <j } [L_i, L _j ] x_i x_j.$$  Since $q=(q_1, \dots, q_n)$
is in the bijective form, i.e.,
$$ q_i = x_i ^2 + \sum _{l<k} \gamma ^{(i)} _{lk} x_l x_k$$
for each $i$, we have
$$ \sum _i K_i q_i = \sum _i K_i x_i ^2 + \sum _i \sum _{l<k} K_i
\gamma ^{(i)} _{lk} x_l x_k.$$ Comparing the coefficient of $x_i
^2$, we get $K_i = L_i ^2 +L _i$ for each $i=1,\dots, n$. This
completes the proof.
\end{proof}

An immediate consequence of Proposition~\ref{pro:specialform} is
the following:

\begin{pro}
\label{pro:commdiagram} Let $E(Q): 0 \to V \to G(Q) \to W \to 0$
be a Bockstein closed 2-power exact sequence corresponding to the
quadratic map $Q: W \to V$. Suppose that we have chosen a
bijective basis for $V$ and $W$ as in Lemma~\ref{lem:nicebasis},
and let  $ L = \sum _i L_i x_i$ be the binding matrix with respect
to this basis. Then, the following diagram commutes
$$
\begin{CD}
W @>\rho>> \gl (V) \\
@VVQV @VV{Q_{\gl}}V \\
V @>\xi>> \gl (V) \\
\end{CD}
$$
where $\rho $ and $\xi$ are linear transformations such that $
\rho (w_i ) = L_i $ and $\xi (v_i) = L_i ^2 + L_i $ for all
$i=1,\dots, n$.
\end{pro}

\begin{proof} We will show that the quadratic maps $\xi \circ Q $ and
$Q_{\gl} \circ \rho $ are equal by showing that the corresponding
cohomology classes  $\xi _* (q) $ and $ \rho ^* (q_{\gl} ) $ are
equal where
$$
\begin{CD}
H^2 (W, V ) @>\xi _*>> H^2 ( W, \gl (V) )
@<\rho ^*<< H^2 (\gl (V), \gl (V) ) \\
\end{CD}
$$
Note that $\xi _* (q)= \sum _i (L_i ^2 + L_i )q_i $. On the other
hand $\rho ^* (q_{\gl } )$ is equal to $ \rho ^* ( \beta (\MX) +
\MX ^2 ) = \beta (L) + L^2 $ (see Example \ref{ex:glagain}). By
Proposition~\ref{pro:specialform} these two are equal.
\end{proof}

The commuting diagram given in Proposition \ref{pro:commdiagram}
can be thought of as a representation of a quadratic map. In terms
of group extensions we have:

\begin{cor}
\label{cor:morpofextensions} The following diagram commutes:
$$
\begin{CD}
E(Q): 0 @>>> V @>>> G(Q) @>>> W @>>> 0 \\
@. @VV{\xi}V @VV{\varphi}V @VV{\rho}V @. \\
E(Q_{\gl}) : 0 @>>> \gl(V) @>>>
K_n(\ZZ / 8 ) @>>> \gl (V)  @>>> 0 \\
\end{CD}
$$
where $\xi$ and $\rho$ are as above.
\end{cor}

We will usually refer to such a diagram as a morphism of group
extensions. We can think of this morphism as a representation of
$E(Q)$ into a $\gl$-induced extension. Whenever there is such a
morphism, we can define its kernel as the extension
$$ 0 \to \ker \xi \to \ker \varphi \to \ker \rho \to 0 $$
whose associated quadratic map is given by $ Q |_{\ker \rho} :
\ker \rho \to \ker \xi .$ Note that the representation given by
the binding operator in general has a non-trivial kernel. When
this happens one can express the original extension as an
extension of group extensions with smaller rank. The following is
an example of this situation:

\begin{ex}
\label{ex:noninjectivemorph} Let $E(Q)$ be the extension with
extension class
$$q=(x^2 +yz,\ y^2+xz,\ z^2).$$ This extension was considered
earlier in Example~\ref{ex:Rusincounterex}, and using the
calculations there, it is easy to see that  the image of both
$\rho$ and $\xi$ have rank $2$, so the image of $\varphi$ is an
extension of the form
$$E(Q'): 0 \to (\ZZ / 2 )^2 \to {\rm im}\ \varphi \to (\ZZ / 2 )^2
\to 0.$$ Further calculations show that the extension class for
$E(Q')$ is $q'=(s^2 +st,\ t^2 )$. In this case, the kernel is the
extension
$$E(Q''): 0 \to  \ZZ /2 \to \ker \varphi \to \ZZ /2  \to 0 $$
where $q''=u^2$.
\end{ex}

\begin{rem} Note that the situation given in Example
\ref{ex:noninjectivemorph} can be best described by saying that we
have a ``central extension" of group extensions of the form
$$ 0 \to E(Q'') \to E(Q) \to E(Q') \to 0.$$
It is possible to develop an extension theory for group extensions
(or for quadratic maps) together with an appropriate
representation theory and make sense of this. We leave this to
another paper which we plan to write as a sequel to this paper.
\end{rem}

Another observation we can make about Example
\ref{ex:noninjectivemorph} is that the quadratic map $q$ is
triangulable since both $q''$ and $q'$ are triangulable. In fact,
the triangulation will come from a basis $\{ v_1, v_2, v_3 \}$
such that $v_1$ generates the kernel of $\xi$. One may ask if in
general the morphism of extensions obtained from the binding
operator gives us a triangulation (possibly by applying the
procedure repeatedly).

We should also mention here that we do not have any examples of
non-triangulable Bockstein closed $2$-power exact extensions.
There is not enough evidence to claim that every Bockstein closed
$2$-power exact extension is triangulable, but it is certainly
tempting to ask if it is true. The following shows that to find a
non-triangulable Bockstein closed extension counterexample one
needs at least $\dim V \geq 3$.

\begin{pro}
\label{pro:n=2 case} Let $E: 0 \to V \to G \to W \to 0$ be a
Bockstein closed extension. If $\ \dim V=2$, then $E$ is
diagonalizable.
\end{pro}

\begin{proof}
Let $q=(q_1, q_2 )$ be the extension class for $E$. We can assume
the extension is Frattini, because otherwise the result holds for
trivial reasons. Since $E$ is Bockstein closed, there is an
$$L=\left[\begin{matrix} f & g \cr h & k  \cr
\end{matrix}\right] \ \  {\rm such \ that}  \ \ \ \
\left[\begin{matrix}\beta (q_1) \cr \beta
(q_2)\end{matrix}\right]=\left[\begin{matrix} f & g \cr h & k \cr
\end{matrix}\right]\left[\begin{matrix} q_1 \cr q_2
\end{matrix}\right].$$
Applying the Bockstein operator again, we get
$$[\beta (L)+L^2]q=\left[\begin{matrix} gh & g(g+f+k)
\cr h(h+f+k) & gh \cr
\end{matrix}\right]\left[\begin{matrix} q_1 \cr q_2
\end{matrix}\right]=0.$$
If $g=0$, then the above equation gives $h(h+f+k)q_1=0$, which
implies $h(h+f+k)$ by Frattini condition. So, $\beta (L) +L^2=0$,
and hence $L$ is diagonalizable. Similarly,  one can show that $L$
is diagonalizable when $h=0$ as well.  So, assume both $g$ and $h$
are non-zero. Then, we get
$$\left[\begin{matrix} h & g+f+k \cr h+f+k & g \cr
\end{matrix}\right]\left[\begin{matrix} q_1 \cr q_2
\end{matrix}\right]=0.$$
This gives, $(f+k)(q_1+q_2)=0$, so we get $f=k$ by Frattini
condition. Setting $f+k=0$ in the above matrix equation, we get
$hq_1+gq_2=0$. Using this we get
$$\left[\begin{matrix}\beta (q_1) \cr \beta
(q_2)\end{matrix}\right]=\left[\begin{matrix} f+h & 0 \cr 0 & g+k
\cr
\end{matrix}\right]\left[\begin{matrix} q_1 \cr q_2
\end{matrix}\right].$$
Hence, $E$ is diagonalizable.
\end{proof}

\section{Strongly Bockstein Closed Exact Sequences}
\label{sect:StronglyB-closed}

The main examples of Bockstein closed quadratic maps are $\gl
$-induced extensions $$E(Q): 0 \to W \to G(Q) \to W \to 0$$ where
$W $ is a square closed subspace of $\gl _n (\FF _2 )$ and  the
quadratic map $Q: W \to W$ is given by $Q(\MA)=\MA ^2 +\MA$. In
this case the associated bilinear map $B(\MA, \MB ) =\MA \MB + \MB
\MA$ satisfies the equation
$$B(Q(\MA), \MB) =B(\MA
,\MB)+B(B(\MA, \MB), \MA)$$ so the equation in
Theorem~\ref{thm:B-Ptheorem} holds with $P=B$. These extension are
Bockstein closed in a stronger sense.

\begin{defn}
\label{defn:StronglyB-closed} Let $Q: W \to W $ be a quadratic map
and $B:W \times W \to W $ be the associated bilinear map. We say
$Q$ is {\bf strongly Bockstein closed} if the equation
$$B(Q(x), y ) =B(x,y) +B(B(x,y),x)$$
holds for all $x,y \in W$. Equivalently, a quadratic map is
strongly Bockstein closed if the equation $\beta (q)= L q$ holds
with $L=Ad_B$ where $Ad_B: W \to \End (W)$ is the adjoint operator
defined as $Ad_B(x)(y)=B(x,y)$.
\end{defn}

There is a close connection between strongly Bockstein closed
extensions and $2$-restricted Lie algebras. To explain this
connection, we first recall the definition of $2$-restricted Lie
algebras.

\begin{defn}
\label{defn:resLiealg} Let $W$ be a vector space over $\FF _2$. We
say $\mathfrak{g}= (W, [\cdot , \cdot], (\cdot ) ^{[2]})$ is a
$2$-restricted Lie algebra if the following holds for all $x,y \in
W$,\\
$(i) \ (x+y)^{[2]} = x^{[2]}+y^{[2]}+[x,y], $\\
$(ii) \ [x^{[2]}, y ]=[[x,y],x].$
\end{defn}

Let \ $Q:W \to W$ be a strongly Bockstein closed quadratic map,
and let $B:W \times W \to W$ be the bilinear map associated to
$Q$. If we take $B$ as a bracket and define the $2$-power map on
$W$ by the formula $w^{[2]}=Q(w)+w$, then $W$ becomes a restricted
Lie algebra together with this bracket and $2$-power map.
Conversely, if $\mathfrak{g}= (W, [\cdot , \cdot], (\cdot )
^{[2]})$ is a $2$-restricted Lie algebra, then the quadratic map
defined by $Q(w)=w^{[2]}+w$ is a strongly Bockstein closed
quadratic map. We conclude the following:

\begin{pro} There is a one-to-one correspondence between strongly
Bockstein closed extensions and  $2$-restricted Lie algebras
$\mathfrak{g}= (W, [\cdot , \cdot], (\cdot ) ^{[2]})$. The
correspondence is given by the formulas $Q(w)=w^{[2]}+w$ and
$B(w,w')=[w,w']$ where $B$ is the bilinear map associated to $Q$.
\end{pro}

Recall that a representation of a $2$-restricted Lie algebra is a
linear transformation $$ \rho : W \to \gl (V) $$ which commutes
with brackets and $2$-power maps where the $2$-power map on $\gl
(V)$ is given by $A \to A^2$. Note that given a representation
$\rho : W \to \gl (V)$ of $2$-restricted Lie algebras, we have a
commuting diagram of the form
$$
\begin{CD}
W @>\rho>> \gl (V) \\
@VVQV @VVQ_{gl}V \\
W @>\rho>> \gl (V) \\
\end{CD}
$$
where $Q$ is the associated quadratic map to the $2$-power map on
$W$, and $Q_{\gl}$ is the quadratic map $A \to A^2 +A$ as before.

Since every $2$-restricted Lie algebra $\mathfrak{g}$ has a
faithful restricted Lie algebra representation $\rho : W \to \gl
_k (\FF _2 )$ for some $k$ (see for example page 192 of
\cite{Jacobson}), every strongly Bockstein closed quadratic map
has a commuting diagram of the above form with injective
horizontal maps. From this, we conclude the following:

\begin{pro}
\label{pro:allglninduced} If \ $Q:W \to W$ is a strongly Bockstein
closed quadratic map, then it is a $\gl$-induced quadratic map.
\end{pro}

In the rest of the section we consider strongly Bockstein closed
extensions which are $2$-power exact. It turns out that these are
the true analogues of Bockstein closed $p$-power exact extensions
with $p>2$. We now explain this analogy:

Let $E:0 \to V \to G \to W \to 0$ be a $p$-power exact extension
where $V$ and $W$ are $\FF_p$-vector spaces with $p>2$. Recall
that $E$ is $p$-power exact means that $\dim W=\dim V$ and that
$E$ is both Frattini and effective. The cohomology of these
extensions has been studied by Browder and Pakianathan in
\cite{BrowPak}. One of the main ingredients for their analysis is
the fact that the $p$-power map gives an isomorphism $W \to V$
under the assumption of $p$-power exactness. Using the inverse of
the $p$-power map, they turned the bilinear commutator map
$[\cdot, \cdot ] : W \times W \to V$ into a bracket $[\cdot, \cdot
]: W \times W \to W$. Then they showed that the extension is
Bockstein closed if and only if this bracket is a Lie bracket.

For $p=2$, the general situation is much more complicated, but the
case of strongly Bockstein closed extensions is, in fact, very
similar to the $p>2$ case. As before let $$E(Q): 0 \to V \to G(Q)
\to W \to 0$$ denote a $2$-power exact sequence corresponding to
the quadratic map $Q: W \to V$. By Lemma~\ref{lem:nicebasis}, we
can choose a bijective basis $\{w_1, \dots, w_n\}$ and $\{
v_1,\dots , v_n \}$ such that $Q(w_i)=v_i$ for all $i= 1, \dots ,
n $. Let $\{x _1 , \dots , x_n \}$ be the basis dual to $\{w_1,
\dots, w_n\}$, and let $\phi : V \to W $ be the linear map which
takes $v_i $ to $w_i$ for all $i=1, \dots, n$. We can replace $Q$
with a quadratic map $\phi \circ Q: W \to W $ and bilinear map $B$
with $\phi \circ B : W \times W \to W .$ We will still denote them
with $Q$ and $B$ as before. With respect to the basis $\{w_1,
\dots, w_n\}$, we have the following matrix representation
\begin{equation}
\label{eqn:main}
q = \beta(x)+ Ax
\end{equation}
where $x$ is the column matrix with $i$-th entry equal to $x_i$,
and $A$ is a matrix whose entries are linear polynomials in
$x_i$'s. In fact, we can describe $A$ in terms of bilinear form
$B$ as follows:

Let $A=\sum _{i=1} ^n A_i x_i $, then $A_i (j,k) =B_j (w_i, w_k )$
if $i< k$ and zero otherwise. Similarly, we can define $A^{\perp}$
as the matrix $A^{\perp}=\sum _{i=1} ^n A ^{\perp} _i x_i $ where
$A^{\perp} _i (j,k) =B_j (w_i, w_k )$ if $i>k$ and zero otherwise.
Note that $A+A^{\perp}$ is the matrix for the adjoint map $Ad_B : W
\to \Hom (W, W)$ defined by $Ad_B (w) (w')=B(w,w')$.

Applying Bockstein to the Equation \ref{eqn:main}, we get
$$ \beta (q)=\beta (A) x+A \beta (x)=\beta (A)x+A^2 x+ Aq$$
We can write $\beta (A)x= A^{\perp} \beta(x) =A^{\perp}q+
A^{\perp} Ax$. So we have
$$ \beta (q)= (Ad_B)Ax + (Ad_B)q.$$

For $p>2$, the above equation holds only when  $(Ad _B)Ax=0$ and
$\beta (q)=(Ad _B )q$. This has two consequences: One is that for
$p>2$ an extension is Bockstein closed only if the $Ad_B$ is a
binding operator, i.e., Bockstein closeness is the same as strongly
Bockstein closeness.  The second consequence is that an extension is
Bockstein closed only if $Ad_B$ satisfies a specific equation  $(Ad
_B)Ax=0$. So it is much easier to decide on Bockstein closeness for
$p>2$. Moreover, it turns out that in this case the equation $(Ad
_B)Ax=0$ is equivalent to the fact that the bilinear map $B : W
\times W \to W$ is a Lie bracket, giving the result of Browder and
Pakianathan  \cite{BrowPak} that the extension is Bockstein closed
if and only if the associated bracket $B$ is a Lie bracket.

In the case $p=2$,  we do not have such a strong conclusion. In
general, one can have Bockstein closed extensions without the
equation $(Ad _B)Ax=0$ holding. However, in the case of strongly
Bockstein closed extensions, one has $\beta (q)=Ad _B q $, forcing
$(Ad_B)Ax=0$ to hold. For $p=2$, the equation $(Ad_B)Ax=0$ is
equivalent to the condition that $B$ is a Lie bracket for a
$2$-restricted Lie algebra with $2$-power map given by
$w^{[2]}=Q(w)+w$. So, as in the case of $p>2$, we have a
connection with Lie algebras, but instead of with the usual Lie
algebras, with restricted Lie algebras. Thus:

\begin{lem}
\label{lem:correspondence} There is a 1-1 correspondence between
strongly Bockstein closed $2$-power exact quadratic maps $Q: W \to
W$, and $2$-restricted Lie algebras $\mathfrak{g}=(W,\ [\cdot,
\cdot ],\ (\cdot )^{[2]})$ satisfying
the following properties: \\
(i) The elements of the form $w^{[2]}+w$ generate $W$, \\
(ii) for all $w\in W$, $w^{[2]}=w$ implies $w=0$. \\
The correspondence is given by the equation $Q(w)=w^{[2]}+w$.
\end{lem}

If a $2$-restricted lie algebra satisfies the properties listed in
the above lemma, we call this Lie algebra a {\bf $2$-power exact
restricted Lie algebra}. It is easy to give an example of such Lie
algebras.

\begin{ex}
\label{ex:mainex} We have seen earlier that the vector space of $n
\times n$ strictly upper triangular matrices, denoted by
$\mathfrak{u} _n$, forms a $2$-restricted lie algebra such that
the quadratic map defined by $Q(\MA)=\MA^2 +\MA$ is effective and
Frattini. So, $\mathfrak{u} _n$ is a $2$-power exact restricted
Lie algebra.
\end{ex}

We would like to remark that $2$-power exact restricted Lie
algebras are quite special. In particular we have the following
result:

\begin{pro}
\label{pro:upper} Let $\mathfrak{g}$ be a $2$-power exact
restricted Lie algebra. If $\ \mathfrak {g}$ has a filtration with
ideals
$$0=\mathfrak{g}_0 \subseteq \mathfrak{g}_1 \subseteq \cdots
\subseteq \mathfrak{g}_k=\mathfrak{g}$$ such that $\dim
(\mathfrak{g}_i/\mathfrak{g}_{i-1} )=1$, then $\mathfrak{g}$ is
nilpotent.
\end{pro}

\begin{proof}
Let $\mathfrak{g}_1=\la x_0 \ra $. There is a linear map $\alpha :
\mathfrak{g} \to \FF _2 $ such that $[x_0, x]=\alpha (x) x_0$ for
every $x \in \mathfrak{g}$. If $x=w^{[2]}+w$ for some $w$, then $$
[x_0, x]=[x_0, w^{[2]}+w]= [[x_0, w], w]+[x_0,w]=[\alpha
(w)^2+\alpha (w)] x_0=0.$$ Since the elements of the form
$w^{[2]}+w$ generate $W$, we can conclude that $x_0$ is central.
So, the result follows from the following lemma.
\end{proof}

\begin{lem}
\label{lem:inductionstep} Let $\mathfrak{g}$ be a $2$-power exact
restricted Lie algebra, and let $Z(\mathfrak{g})$ be the center of
$\mathfrak{g}$. Then, $\mathfrak{g}/ Z(\mathfrak{g})$ is also a
$2$-power exact restricted Lie algebra with $2$-power map defined
by $x+Z(\mathfrak{g}) \to x^{[2]}+ Z(\mathfrak{g})$.
\end{lem}

\begin{proof}
It is easy to see that $\mathfrak{g}/ Z(\mathfrak{g})$ satisfies
the Frattini condition. Now we will show that it is effective.
Assume $x+x^{[2]}=z$ for some $z \in Z(\mathfrak{g})$. Since the
map defined by $x \to x^{[2]}+x$ is linear on $Z(\mathfrak{g})$,
it must be an isomorphism by the effectiveness condition. So, we
can find a $y \in Z(\mathfrak{g})$ such that $y^{[2]}+y=z$. Then,
we have
$$(x+y)^{[2]}=x^{[2]}+y^{[2]}=(x+z)+(y+z)=x+y.$$ By the effectiveness
of $\mathfrak{g}$, we must have $x+y=0$. So, $x \in
Z(\mathfrak{g})$ giving that $\mathfrak{g}/ Z(\mathfrak{g})$ is
$2$-power exact.
\end{proof}

We do not know if all the $2$-power exact restricted Lie algebras
are nilpotent. It would be interesting to see whether the general
theory of restricted Lie algebras applies to answer this question.
We leave this as an open problem.

{\it Acknowledgments:} We thank the referee for his comments on the
paper.

\bigskip

\noindent
Dept. of Mathematics \\
University of Rochester, \\
Rochester, NY 14627 U.S.A. \\
E-mail address: jonpak@math.rochester.edu \\

\bigskip

\noindent
Dept. of Mathematics\\
Bilkent University,\\
Ankara, 06800, Turkey. \\
E-mail address: yalcine@fen.bilkent.edu.tr \\

\end{document}